\documentclass[12pt]{amsart}
\usepackage{amssymb,latexsym, amscd, a4wide}
\usepackage[all]{xy}
\usepackage{graphicx}
\usepackage{mathrsfs}
\usepackage{amsmath}
\usepackage{pb-diagram}
\usepackage{bbold}
\usepackage{color}

\newtheorem{thm}{Theorem}[section]
\newtheorem{lem}[thm]{Lemma}
\newtheorem{prop}[thm]{Proposition}

\newtheorem{defi}[thm]{Definition}
\newtheorem{rem}[thm]{Remark}

\newtheorem*{thm*}{Theorem}
\newtheorem*{prop*}{Proposition}

 \numberwithin{equation}{section}
 \newcommand{\nc}{\newcommand}
\nc{\ssn}{\subsection{}} \nc{\sssn}{\subsubsection{}}

\newcommand{\dem }{\noindent\textbf{Proof}. }
\newcommand{\findem }{\hfill $\Box$ \vskip0.5cm }

\numberwithin{equation}{section}
\newcommand{\Ind}{\operatorname{Ind}}

\nc{\oR}{\ol{R}}
\newcommand{\ol}{\ensuremath{\overline}}

\nc{\id}{\textrm{id}}

\newcommand{\N}{\mathcal{N}}

\newcommand{\Ext}{\operatorname{Ext}}

\newcommand{\ad}{\hbox{\ensuremath{\operatorname{ad}}}}

\newcommand{\Z}{\ensuremath{\mathbb{Z}}}

\renewcommand{\mod}{\ensuremath{\operatorname{mod}}}

\newcommand{\C}{\ensuremath{\mathbb{C}}}

\newcommand{\g}{\ensuremath{\mathfrak{g}}}
\newcommand{\h}{\ensuremath{\mathfrak{h}}}

\newcommand{\Sym}{\operatorname{Sym}}

\newcommand{\mc}[1]{\mathcal{#1}} 
\newcommand{\mb}[1]{\mathbb{#1}} 
\newcommand{\mf}[1]{\mathfrak{#1}} 

\begin{document}
\title[Crystal bases for reduced imaginary Verma modules]{Crystal bases for reduced imaginary Verma modules of untwisted quantum affine algebras}
\author{Juan Camilo Arias}
\author{Vyacheslav Futorny}
\author{Kailash C. Misra}

\address{Institute of Mathematics and Statistics, University of S\~ao Paulo, S\~ao Paulo, Brazil}
\email{jcarias@ime.usp.br }
\address{Shenzhen International Center for Mathematics, Southern University of Science and Technology, China, and Institute of Mathematics and Statistics,
 University of S\~ao Paulo,
 S\~ao Paulo, Brazil}
 \email{vfutorny@gmail.com}
 \address{Department of Mathematics, 
 North Carolina State University, 
 Raleigh, NC, USA}
 \email{misra@ncsu.edu}


\maketitle
\begin{center}
{\it To the memory of Georgia Benkart and Ben Cox}
\end{center}

\begin{abstract}  
We consider reduced imaginary Verma modules for the untwisted quantum affine algebras $U_q(\hat{\g})$ 
 and define a crystal-like base which we call imaginary crystal base 
 using  the Kashiwara algebra $\mathcal K_q$ constructed in earlier work by Ben Cox and two of the authors. We prove the existence of the imaginary crystal base 
 for any object in a suitable category $\mc{O}^q_{red,im}$ containing the reduced imaginary Verma modules for $U_q(\hat{\g})$. 
\end{abstract}

\section{Introduction}
Let $\hat{\g}$ be an untwisted affine Lie algebra with Cartan subalgebra $\hat{\h}$ and generalized Cartan matrix $A=(a_{ij})_{0\leq i,j\leq N}$. Let $\g$ be the associated finite dimensional simple Lie algebra with Cartan matrix $(a_{ij})_{1\leq i,j\leq N}$ and Cartan subalgebra $\h$. Let $\{\alpha_0, \alpha_1, \ldots, \alpha_N\}$ be the set of simple roots, $\delta$ be the null root and $\Delta$ be the set of roots for $\hat{\g}$ with respect to $\hat{\h}$. A partition $\Delta = S \cup -S$  is a closed partition if $\alpha , \beta \in S$ and 
$\alpha + \beta \in \Delta$ implies $\alpha + \beta \in S$. The classification of closed partitions of the affine root system was obtained by Jakobsen and Kac \cite{JK89}, and independently by Futorny \cite{Fut90, Fut92}. The usual partition of the set of roots $\Delta = \Delta_+ \cup \Delta_-$ to the set of positive and negative roots is called the standard partition. Corresponding to this partition we have a standard Borel subalgebra of $\hat{\g}$ from which we induce the standard Verma modules. Consider the nonstandard closed partition
$\Delta = S \cup -S$, where $S = \{ \alpha + n\delta | \alpha\in \Delta_{0,+}, n\in \Z   \} \cup \{ k\delta | k>0   \}$, $\Delta_{0,+}$ being the positive roots for $\g$. Corresponding to this partition we obtain an inequivalent Borel subalgebra of $\hat{\g}$ which induce nonstandard Verma modules $M(\lambda)$ called imaginary Verma modules. Unlike the standard Verma modules, the imaginary Verma modules contains both finite and infinite dimensional weight spaces. As shown in \cite{CFKM, FutGM} these imaginary Verma modules can be $q$-deformed to quantum imaginary Verma modules $M_q(\lambda)$ for the quantum affine algebras $U_q(\hat{\g})$ preserving both finite and infinite dimensional weight spaces for generic $q$. 
\vskip 5pt
The theory of crystal bases for integrable representations of $U_q(\hat{\g})$ was developed independently by Kashiwara \cite{Kas90, Kas91} and Lusztig \cite{L03}. Two of the authors (KCM and VF) jointly with Ben Cox initiated the investigation of the existence of crystal-like bases for quantum imaginary Verma modules $M_q(\lambda)$ for the quantum affine algebras $U_q(\hat{\g})$ in 2008 at a conference in Banff. Following the framework in \cite{Kas91} they constructed an analog of Kashiwara algebra $\mathcal{K}_q$ for $M_q(\lambda)$ in the case $\hat{\g} = \hat{sl}(2)$ by introducing certain Kashiwara type operators in \cite{CFM} and proved that certain quotient $\mathcal{N}_q^-$ of $U_q(\hat{\g})$ is a simple $\mathcal{K}_q$-module. This result was extended to all quantum affine algebras $U_q(\hat{\g})$ of $ADE$ types in \cite{CFM01}. In \cite{CFM03} a category $\mc{O}^q_{red,im}$ of $U_q(\hat{sl}(2))$-modules were introduced and it was shown that any module in this category is a simple reduced quantum imaginary Verma module $\tilde{M}_q(\lambda)$, certain quotient of $M_q(\lambda)$  or direct sum of these modules. The existence of imaginary crystal base for any module in the category $\mc{O}^q_{red,im}$ was shown in \cite{CFM03, CFM04}.
\vskip 5pt
In this paper, by appropriate modifications we first extend the results in \cite{CFM01} to all untwisted quantum affine algebras $U_q(\hat{\g})$. Using the results in \cite{AFO}, we construct the category $\mc{O}^q_{red,im}$ of $U_q(\hat{\g})$-modules and show that any module in this category is a simple reduced quantum imaginary Verma module $\tilde{M}_q(\lambda)$ or direct sum of these modules. Moreover,  we show that any module in the category $\mc{O}^q_{red,im}$  admits an imaginary crystal basis.

 From now on imaginary Verma modules will mean the quantum imaginary Verma modules $M_q(\lambda)$ unless there is any confusion. 

\vskip 5pt

This paper is organized as follows. In Sections $2$ and $3$, we define and set the notations for affine and quantum affine algebras. In Section $4$, we recall the definitions of imaginary Verma modules $M_q(\lambda)$ and reduced imaginary Verma modules $\tilde{M}_q(\lambda)$ for any weight $\lambda \in P$ (where $P$ is the weight lattice) and state some of the necessary and sufficient conditions for their irreducibility. In Section $5$, we recall the definitions of $\Omega$-operators and the Kashiwara algebra $\mathcal{K}_q$ from \cite{CFM01} for any untwisted quantum affine algebra $U_q(\hat{\g})$. We show that the unique non-degenerate symmetric form on the quotient $\mathcal{N}_q^-$ of $U_q(\hat{\g})$ defined in \cite{CFM01} holds for any untwisted affine algebra $\hat{\g}$. In Section $6$, we introduce a product between monomials in $\mathcal{N}_q^-$ which we call ``twisted concatenation product", define Kashiwara type operators on $\tilde{M}_q(\lambda)$ and prove some of their relations. In Section $7$ we introduce a bilinear form among ordered monomials of  
$\mathcal{N}_q^-$ and prove that it satisfies certain orthonormality condition modulo $q^2\mathbb{Z}[q]$ which plays an important role in the construction of the imaginary crystal basis. In Section $8$ we prove the existence of imaginary crystal basis for a simple reduced imaginary Verma module $\tilde{M}_q(\lambda)$. Finally, in the last section we define the category $\mc{O}^q_{red,im}$ and show that any module in this category is either a simple reduced imaginary Verma module or a direct sum of these modules. We prove that any module in $\mc{O}^q_{red,im}$ has an imaginary crystal basis.
\vskip 5pt

We dedicate this paper to Georgia Benkart and Ben Cox whom we miss dearly. As mentioned above Ben has been an integral part for the success of this project in the earlier stages. He has been a friend and collaborator of KCM and VF for many years. Georgia Benkart was a research mentor and collaborator of KCM. Both KCM and VF are indebted for her support and encouragements over the years. 
\vskip 5pt

\section{Affine algebras} 

Let $I=\{ 0, 1, \ldots, N\}$ and $A=(a_{ij})_{0\leq i,j\leq N}$  be a generalized affine Cartan matrix for an untwisted affine Kac-Moody Lie algebra $\hat{\g}$. Let $D=diag(d_0, d_1, \ldots, d_N)$ be a diagonal matrix with relatively prime integer entries such that $DA$ is symmetric. The numbers $d_0, \ldots, d_N$ are given as follows: Types ADE: $d_i=1$ for all $i\in I$, type B: $d_i=2$ for $i\in I\setminus\{N\}$ and $d_N=1$, type C: $d_0=d_N=2$ and $d_i=1$ for $i\neq 0,N$, Type F: $d_0=d_1=d_2=2$, $d_3=d_4=1$ and type G: $d_0=d_1=3$ and $d_2=1$.\\ 

The Chevalley-Serre presentation of $\hat{\g}$ is given by generators $e_i, f_i, h_i$ for $0\leq i \leq N$ and $d$ subject to the defining relations:

$$ [h_i,h_j]=0 \quad [d,h_i]=0 \quad [h_i,e_j]=a_{ij}e_j \quad [h_i,f_j]=-a_{ij}f_j $$
$$[e_i,f_j]=\delta_{i,j}h_i \quad [d,e_i]=\delta_{0,i}e_i \quad [d,f_i]=-\delta_{0,i}f_i$$
$$(\ad e_i)^{1-a_{ij}}(e_j)=0 \quad (\ad f_i)^{1-a_{ij}}(f_j)=0.$$

Then the abelian subalgebra $\hat{\h}= {\rm span}\{h_0, \ldots, h_N,d\}$ is the Cartan subalgebra of $\hat{\g}$. Let $\Delta$  be the root system of $\hat{\g}$ with simple roots $\Pi=\{\alpha_0, \alpha_1, \ldots, \alpha_N\}$ and let $\delta=\alpha_0+\theta$ be the null root where $\theta$ is the longest root of the underlying simple Lie algebra $\g$. Recall that $\Delta = \Delta^{re} \cup \Delta^{im}$, where $\Delta^{re}$ and $\Delta^{im}$ denotes the real and imaginary sets of roots. Let $Q$, $P$, $\check{Q}$ and $\check{P}$ denote the root lattice, the weight lattice, the coroot lattice and the coweight lattice respectively.
The standard non-degenerate symmetric bilinear form $(.|.)$ on $\hat{\h}$ satisfies $(\alpha_i|\alpha_j)=d_ia_{ij}$, $(\delta|\alpha_i)=(\delta|\delta)=0$ for all $i,j\in I$.\\

Let $I_0=\{ 1,2, \ldots, N\}$ and let $\g$ be the simple finite dimensional Lie algebra with Cartan matrix  $(a_{ij})_{1\leq i,j \leq N}$. Let $\Delta_0$ be the root system of $\g$ and $\Pi_0=\{\alpha_1, \ldots, \alpha_n\}$ be the set of simple roots. As a matter of notation, any subscript $0$ will refers to the same sets for the simple Lie algebra $\g$.\\

The loop space realization of $\hat{\g}$ is given by $\hat{\g} = \g\otimes \C[t,t^{-1}]\oplus \C c \oplus \C d$, where $c$ is a central element, $d$ is a degree derivation such that $[d,x\otimes t^n]=nx\otimes t^n$ for any $x\in \g$ and $n\in \Z$, and 
$$[x\otimes t^n, y\otimes t^m] = [x,y]\otimes t^{n+m} + \delta_{n,-m}n(x|y)c,$$ for all $x,y\in\g$, $n,m\in\Z$. Here $(-|-)$ is a symmetric invariant bilinear form on $\hat{\g}$.\\

Recall that a subset $S$ of $\Delta$ is said to be closed if whenever $\alpha, \beta \in S$ and $\alpha + \beta \in \Delta$ then $\alpha + \beta \in S$. We also say that $S$ is a closed partition if $S$ is closed, $\Delta = S \cup -S $ and $S\cap -S = \emptyset$. We are going to consider the closed partition, inequivalent to the standard partition into positive and negative roots, given by  $S = \{ \alpha + n\delta | \alpha\in \Delta_{0,+}, n\in \Z   \} \cup \{ k\delta | k>0   \}$.

\section{Quantum affine algebras}

Let $U_q(\hat{\g})$ be the quantum affine algebra (see \cite{CP}, \cite{L01}), i.e., the associative and unital $\C(q^{1/2})$-algebra with generators $E_i, F_i, K_{\alpha}, \gamma^{\pm1/2}, D^{\pm 1}$ for $0\leq i \leq N$, $\alpha\in Q$ and defining relations:

$$ DD^{-1} = D^{-1}D = K_\alpha K_{-\alpha} = K_{-\alpha} K_\alpha = \gamma^{1/2}\gamma^{-1/2}=1 $$
$$ [\gamma^{\pm1/2}, U_q(\hat{\g})]=[D,K_{\pm \alpha}]=[K_\alpha,K_\beta]=0$$
$$ (\gamma^{\pm 1/2})^2 = K_{\pm \delta} $$
$$ E_iF_j-F_jE_i = \delta_{ij}\frac{K_i - K_i^{-1}}{q_i+q_i^{-1}} $$
$$ K_{\alpha}E_iK_{-\alpha}= q^{(\alpha|\alpha_i)}E_i, \quad K_{\alpha}F_iK_{-\alpha}= q^{-(\alpha|\alpha_i)}F_i$$
$$ DE_iD^{-1} = q^{\delta_{i,0}}E_i, \quad DF_iD^{-1} = q^{-\delta_{i,0}}F_i $$
$$ \sum_{s=0}^{1-a_{ij}} (-1)^s E_i^{(1-a_{ij}-s)}E_jE_i^{(s)} =  \sum_{s=0}^{1-a_{ij}} (-1)^s F_i^{(1-a_{ij}-s)}F_jF_i^{(s)}, \quad i\neq j $$

 where $q_i = q^{d_i}$, $[n]_i=\frac{q_i^n - q_i^{-n}}{q_i-q_i^{-1}}$, $[n]_i!=[n]_i[n-1]_i\cdots [2]_i[1]_i$,  $K_i=K_{\alpha_i}$, $E_i^{(s)} = E_i^{s}/[s]_i!$ and $F_i^{(s)} = F_i^{s}/[s]_i!.$\\

There is also another realization due to Drinfeld in the spirit of the loop space realization as follows (see \cite{D01}, \cite{B01}). $U_q(\hat{\g})$ is the associative unital $\C(q^{1/2})$-algebra generated by $x_{ir}^{\pm1}, h_{is}, K_i^{\pm 1}, \gamma^{\pm1/2}, D^{\pm 1}$ for $i\in I_0$, $r,s\in \Z$ and $s\neq 0$ subject to the relations:

$$ D^{\pm 1}D^{\mp 1}  = K_i^{\pm 1}K_i^{\mp 1} = \gamma^{\pm 1/2}\gamma^{\mp 1/2} = 1  $$  
$$ [\gamma^{\pm 1/2}, U_q(\hat{\g})] = [D,K_i^{\pm 1}] = [K_i, K_j] = [K_i, h_{js}] = 0    $$
$$ Dh_{ir}D^{-1} = q^rh_{ir}, \quad \quad  Dx_{ir}^{\pm}D^{-1} = q^rx_{ir}^{\pm}$$
$$ K_ix_{jr}^{\pm}K_i^{-1} = q^{\pm (\alpha_i | \alpha_j )}x_{jr}^{\pm} $$
$$ [h_{ik}, h_{jl}] = \delta_{k,-l}\frac{1}{k}[ka_{ij}]_i\frac{\gamma^k - \gamma^{-k}}{q_j - q_j^{-1}} $$
$$ [h_{ik}, x_{jl}^{\pm}] = \pm \frac{1}{k}[ka_{ij}]_i\gamma^{\mp |k|/2}x_{j,k+l}^{\pm}  $$
$$ x_{i,k+1}^{\pm}x_{jl}^{\pm} - q^{\pm(\alpha_i|\alpha_j)}x_{jl}^{\pm}x_{i,k+1}^{\pm}  =  q^{\pm(\alpha_i|\alpha_j)}x_{ik}^{\pm}x_{j,l+1}^{\pm} - x_{j,l+1}^{\pm}x_{ik}^{\pm} $$
$$ [x_{ik}^{+}, x_{jl}^{-}] = \delta_{ij}\frac{1}{q_i - q_i^{-1}}(\gamma^{(k-l)/2}\psi_{i,k+l} - \gamma^{(l-k)/2}\phi_{i,k+l})  $$

where

$$ \sum_{k=0}^{\infty} \psi_{ik}z^k = K_i\exp\Big((q_i-q_i^{-1})\sum_{l>0}h_{il}z^l\Big) $$

$$ \sum_{k=0}^{\infty} \phi_{i,-k}z^{-k} = K_i^{-1}\exp\Big(-(q_i-q_i^{-1})\sum_{l>0}h_{i,-l}z^{-l}\Big) $$

and for $i\neq j$, 

$$ \Sym_{k_1, \ldots, k_{1-a_{ij}}}\sum_{r=0}^{1-a_{ij}}(-1)^r{1-a_{ij}\brack r}_i x_{ik_1}^{\pm} \cdots x_{ik_r}^{\pm}x_{ij}^{\pm}x_{ik_{r+1}}^{\pm} \cdots x_{ik_{1-a_{ij}}}^{\pm} = 0.$$

If we consider the following generating functions 

$$ \phi_i(u) = \sum_{p\in\Z}\phi_{ip}u^{-p}, \quad  \psi_i(u) = \sum_{p\in\Z}\psi_{ip}u^{-p}, \quad x_i^{\pm}(u) = \sum_{p\in\Z}x_{ip}^{\pm}u^{-p}$$ 

the defining relations becomes:

$$ [\phi_i(u), \phi_j(v)] = [\psi_i(u), \psi_j(v)] = 0$$
$$ \phi_i(u)\psi_j(v)\phi_i(u)^{-1}\psi_j(v)^{-1} = g_{ij}(uv^{-1}\gamma^{\-1})/g_{ij}(uv^{-1}\gamma)  $$
$$ \phi_i(u)x_j^{\pm}(v)\phi_i(u)^{-1} = g_{ij}(uv^{-1}\gamma^{\mp 1/2})^{\pm 1}x_j^{\pm}(v)  $$
$$ \psi_i(u)x_j^{\pm}(v)\psi_i(u)^{-1} = g_{ji}(vu^{-1}\gamma^{\mp 1/2})^{\mp 1}x_j^{\pm}(v)   $$
$$ (u-q^{\pm(\alpha_i|\alpha_j)v})x_i^{\pm}(u)x_j^{\pm}(v) =  (q^{\pm(\alpha_i|\alpha_j)u-v})x_j^{\pm}(v)x_i^{\pm}(u)  $$
$$ [x_i^{+}(u), x_j^{-}(v)] = \delta_{ij}(q_i - q_i^{-1})(\delta(u/v\gamma)\psi_i(v\gamma^{1/2}) - \delta(u\gamma/v)\phi_i(u\gamma^{1/2}))   $$

where $g_{ij}(t)=g_{ij,q}(t)$ is the Taylor expansion at $t=0$ of the function $(q^{(\alpha_i|\alpha_j)}t-1)/(t-q^{(\alpha_i|\alpha_j)})$ and $\delta(z)=\sum_{k\in\Z}z^k $.

\section{Imaginary Verma modules}

For the closed partition $\Delta = S \cup -S$ where $S = \{ \alpha + n\delta | \alpha\in \Delta_{0,+}, n\in \Z   \} \cup \{ k\delta | k>0   \}$, we define the following subalgebras of $U_q(\hat{\g})$. \\

\begin{itemize}
\item $U_q^{+}(S)$ generated by $x_{ik}^{+}, h_{il}$ for $i\in I_0$, $k\in \Z$ and $l>0$.
\item $U_q^{-}(S)$ generated by $x_{ik}^{-}, h_{i,-l}$ for $i\in I_0$, $k\in \Z$ and $l>0$.
\item $U_q^{0}(S)$ generated by $K_i^{\pm 1}, \gamma^{\pm 1/2}, D^{\pm 1}$ for $i\in I_0$.
\end{itemize}

\vspace{3mm}




For $\lambda \in P$, a wight module $V$ of $U_q(\hat{\g})$ is called an $S$-weight module with highest weight $\lambda$ if there is a non zero vector $v\in V$ of weight $\lambda$ such that $u^{+}v=0$ for all $u^{+}\in U_q^+(S)\setminus \C(q^{1/2})$ and $V=U_q(\hat{\g})v$.\\

Let us consider the Borel subalgebra $B_q$ of $U_q(\hat{\g})$, 
which is generated by $U_q^{+}(S)\cup U_q^0(S)$. 
Consider now the one dimensional module $\C(q^{1/2})_{\lambda}$, in which the $B_q$-module $U_q^{+}(S)$ acts trivially and if $\mathbb{1}$ is the generator of the module then $K_i^{\pm 1}\mb{1} = q^{\pm\lambda(h_i)\mb{1}}$, $i\in I_0$, $\gamma^{\pm 1/2}\mb{1}= q^{\pm\lambda(c)/2}\mb{1} $ and $D^{\pm 1}\mb{1}= q^{\pm \lambda(d)}\mb{1}$.\\


We define the {\it imaginary Verma module} $M_q(\lambda)$ of weight $\lambda \in P$ as 

$$ M_q(\lambda) := U_q(\hat{\g})\otimes_{B_q}\C(q^{1/2})_{\lambda}. $$

\vspace{2mm}

We have the following theorem (see \cite{Fut}).

\begin{thm} $M_q(\lambda)$ is irreducible if and only if $\lambda(c)\neq0$. \end{thm}

Let us suppose that $M_q(\lambda)$ is reducible, that means $\lambda(c)=0$ and so $\gamma^{\pm 1/2}$ acts by $1$. Denote by $J^q(\lambda)$ the left ideal of $U_q(\hat{\g})$ generated by $x_{ik}^{+}, h_{il}$ for $i\in I_0, k,l\in\Z, l\neq0$ and $K_i^{\pm 1} -q^{\pm\lambda(h_i)}$, $\gamma^{\pm 1/2} -1$ and $D^{\pm 1}- q^{\pm \lambda(d)}$. Set 

$$ \tilde{M}_q(\lambda) = U_q(\hat{\g})/J^q(\lambda). $$

\vspace{2mm}

It is a quotient of $M_q(\lambda)$ and we call it {\it reduced imaginary Verma module}. We have the following theorem (see proposition 7.1 in \cite{FGM02}).

\begin{thm} Let $\lambda\in P$ such that $\lambda(c)=0$. Then $\tilde{M}_q(\lambda)$ is irreducible if and only if $\lambda(h_i)\neq 0$ for all $i\in I_0$.  \end{thm}

Notice also that 

$$ \tilde{M}_q(\lambda) = \bigoplus_{\substack{i_1, \ldots, i_r \\ k_1, \ldots, k_r}} \C(q^{1/2})x_{i_1k_1}^-\cdots x_{i_rk_r}^- v_{\lambda}$$

\vspace{2mm}

where $v_{\lambda}$ stands for the generator $\mb{1}$. The order in the defining monomial for the reduced imaginary Verma module is going to be important for us. We will say that a monomial $x_{i_1k_1}^-\cdots x_{i_rk_r}^- $ is ordered if and only if $i_1+k_1 \geq i_2 + k_2\geq \cdots \geq i_r+k_r$. \\

\section{$\Omega$-operators and the Kashiwara algebra}

The Kashiwara algebra and the $\Omega$-operators were defined in \cite{CFM01} for types ADE, but a close look at the proof of the results presented in \cite{CFM01} shows that everything there is true for any quantum affine algebra associated with an untwisted affine Lie algebra of any type. In this section we recall the main constructions presented in \cite{CFM01}.\\

Consider the subalgebra $\mc{N}_q^-$ of $U_q(\hat{\g})$ generated by $\gamma^{\pm 1/2}$ and $x_{il}^-$ for $l\in\Z$, $i\in I_0$ and the relations between them. The elements in $\mc{N}_q^-$ are of the form $P^{j_1, \ldots, j_k}_{m_1, \ldots, m_k} = \gamma^{l/2}x_{j_1m_1}^-\cdots x_{j_km_k}^-$, for $m_i,l\in\Z$, $k\geq0$ and $1\leq j_i\leq N$. Such an element is a summand of 

$$ P^{j_1, \ldots, j_k} = P^{j_1, \ldots, j_k}(v_1, \ldots, v_k) := \gamma^{l/2}x_{j_1}^-(v_1)\cdots x_{j_k}^-(v_k) .$$

\vspace{2mm}

We set $\overline{P}^{j_1, \ldots, j_k} = x_{j_1}^-(v_1)\cdots x_{j_k}^-(v_k)$ and $\overline{P}^{j_1, \ldots, j_k}_l = x_{j_1}^-(v_1)\cdots x_{j_{l-1}}^-(v_{j_{l-1}})x_{j_{l+1}}^-(v_{j_l+1})\cdots x_{j_k}^-(v_k)$.\\

We denote

$$ G_{il} = G_{il}^{1/q} = G_{il}^{1/q}(v_{j1}, \ldots, v_{j_l}, v_l) := \delta_{i,j_l} \prod_{m=1}^{l-1} g_{i,j_m,q^{-1}}(v_{j_m}/v_l) $$ 

$$ G_{il}^q = G_{il}^q(v_{j1}, \ldots, v_{j_l}, v_l) := \delta_{i,j_l}\prod_{m=1}^{l-1} g_{i,j_m,q}(v_l/v_{j_m}) $$

\vspace{2mm}

where $G_{i1}= \delta_{i,j_1}$. We define operators $\Omega_{\psi_i}(k), \Omega_{\phi_i}(k):\mc{N}_q^- \to \mc{N}_q^-$ for $k\in\Z$ in terms of the generating functions

$$ \Omega_{\psi_i}(u) = \sum_{l\in\Z}\Omega_{\psi_i}(l)u^{-l}, \quad \quad \Omega_{\phi_i}(u) = \sum_{l\in\Z}\Omega_{\phi_i}(l)u^{-l} $$

\vspace{2mm}

by 

$$ \Omega_{\psi_i}(u)(\overline{P}^{j_1, \ldots, j_k}) = \sum_{l=1}^k G_{il}\overline{P}^{j_1, \ldots, j_k}_l\delta(u/v_l\gamma) $$

$$\Omega_{\phi_i}(u)(\overline{P}^{j_1, \ldots, j_k}) = \sum_{l=1}^k G_{il}^q\overline{P}^{j_1, \ldots, j_k}_l\delta(u\gamma/v_l). $$

\vspace{2mm}

We have that $\Omega_{\psi_i}(u)(1) = \Omega_{\phi_i}(u)(1)=0$.\\

Recall that $x^{-}(v) = \sum_{m}x_{im}^-v^{-m}$. So we can consider left multiplication operators $x_{im}^-:\mc{N}_q^-\to \N_q^-$. There are several identities between the $\Omega$-operators and the $x^-$-operators, we just list a few:

\begin{multline*}
 q^{(\alpha_i|\alpha_j)}\gamma\Omega_{\psi_j}(m)x_{i,n+1}^{-} - \Omega_{\psi_j}(m+1)x_{in}^{-} =\\
 (q^{(\alpha_i | \alpha_j)}\gamma - 1)\delta_{ij}\delta_{m,-n-1} + \gamma x_{i,n+1}^{-}\Omega_{\psi_j}(m)  - q^{(\alpha_1 | \alpha_j)}x_{in}^{-}\Omega_{\psi_j}(m+1). 
 \end{multline*}
 
 \begin{multline*}
  q^{(\alpha_i|\alpha_j)}\Omega_{\phi_j}(m)x_{i,n+1}^{-} - \gamma\Omega_{\phi_j}(m+1)x_{in}^{-} =\\
 (q^{(\alpha_i | \alpha_j)}-\gamma)\delta_{ij}\delta_{m,-n-1} + x_{i,n+1}^{-}\Omega_{\phi_j}(m)  - q^{(\alpha_1 | \alpha_j)}\gamma x_{in}^{-}\Omega_{\phi_j}(m+1). 
 \end{multline*}

\begin{equation}\label{CommOX} \Omega_{\psi_j}(k)x_{im}^{-} = \delta_{ij}\delta_{k,-m}\gamma^{k} + \sum_{r \geq 0} g_{i,j,q^{-1}}(r)x_{i,m+r}^{-}\Omega_{\psi_j}(k-r)\gamma^r.  \end{equation}

\begin{equation}\label{CommOO} \Omega_{\psi_i}(k) \Omega_{\phi_j}(m) = \sum_{r\geq 0} g_{i,j}(r)\gamma^{2r}\Omega_{\phi_j}(r+m)\Omega_{\psi_i}(k-r). \end{equation}

\vspace{2mm}

We define the {\bf Kashiwara algebra} $\mc{K}_q$ as the $\C(q^{1/2})$-algebra with generators $\Omega_{\psi_j}(m), x_{i}^{-}(n), \gamma^{\pm 1/2}$ for $m,n\in \Z$, $1\leq i,j \leq N$ where $\gamma^{\pm 1/2}$ are central and the defining relations are:

\begin{multline*}
 q^{(\alpha_i|\alpha_j)}\gamma\Omega_{\psi_j}(m)x_{i,n+1}^{-} - \Omega_{\psi_j}(m+1)x_{in}^{-} =\\
 (q^{(\alpha_i | \alpha_j)}\gamma - 1)\delta_{ij}\delta_{m,-n-1} + \gamma x_{i,n+1}^{-}\Omega_{\psi_j}(m)  - q^{(\alpha_1 | \alpha_j)}x_{in}^{-}\Omega_{\psi_j}(m+1) 
 \end{multline*}
 
$$ q^{(\alpha_i \mid \alpha_j)}\Omega_{\psi_i}(k+1)\Omega_{\psi_j}(l) - \Omega_{\psi_j}(l)\Omega_{\psi_i}(k+1) = \Omega_{\psi_i}(k)\Omega_{\psi_j}(l+1) - q^{(\alpha_i | \alpha_j)}\Omega_{\psi_j}(l+1)\Omega_{\psi_i}(k) $$

\begin{equation}\label{OrdX} x_{i,k+1}^- x_{jl}^- -  q^{-(\alpha_i | \alpha_j)}x_{jl}^-x_{i,k+1}^- = q^{-(\alpha_i | \alpha_j)}x_{ik}^-x_{j,l+1}^- - x_{j,l+1}^- x_{ik}^- \end{equation}

and $$\gamma^{\pm 1/2}\gamma^{\mp 1/2} = 1. $$

\vspace{2mm}

Some of the properties of the Kashiwara algebra $\mc{K}_q$ proved in \cite{CFM01} that holds for any quantum untwisted affine algebras are summarized in the following proposition.

\begin{prop}\label{bilinform} For a quantum affine algebra associated to any untwisted affine Lie algebra, there exists a unique non-degenerate symmetric form $(-,-)$ defined on $\N_q^-$ satisfying $(x_{ij}^-a,b)=(a,\Omega_{\psi_i}(-j)b)$ and $(1,1)=1$. Moreover, $\mc{N}_q^-$ is a left $\mc{K}_q$-module such that $\N_q \cong \mc{K}_q / \Big( \sum_{i=1}^n\sum_{k\in\Z} \mc{K}_q\Omega_{\psi_i}(k) \Big)$. Furthermore, $\N_q^-$ is simple as left $\mc{K}_q$-module. \end{prop}

\dem The first statement follows from Proposition 6.0.6 and Corollary 7.0.9 in \cite{CFM01}. The second and third statements follows from Lemma 6.0.5, Theorem 7.0.8 in \cite{CFM01}.
\findem

\section{The operators $\tilde{\Omega}$ and $\tilde{x}$}




Recall that the function $g_{ij,q^{-1}}(r)$ for $i,j\in I_0$ is the Taylor expansion at $t=0$ of the function $(q^{(\alpha_i|\alpha_j)}t-1)/(t-q^{(\alpha_i|\alpha_j)})$ and that $(\alpha_i |\alpha_j) = d_ia_{ij}$. Below we list the explicit form of   $g_{ij,q^{-1}}(r)$ when $r=0$ and $r>0$ for all possible values of $(\alpha_i|\alpha_j)$.

$$ g_{ij,q^{-1}}(r)  = \left\{ \begin{array}{lllll} r=0 & { } & { } & r>0 & \\
q^{6} & { } & { } & (q^{12}-1)q^{6(r-1)} = q^{6(r+1)}(1-q^{-12}) &  (\alpha_i|\alpha_j) = 6\\
q^{4} & { } & { } & (q^{8}-1)q^{4(r-1)} = q^{4(r+1)}(1-q^{-8}) &  (\alpha_i |\alpha_j) = 4\\ 
q^{2} & { }  & { } & (q^{4}-1)q^{2(r-1)} = q^{2(r+1)}(1-q^{-4}) &  (\alpha_i |\alpha_j) = 2\\ 
1 & { } & { } & 0 &  (\alpha_i |\alpha_j) = 0\\ 
q^{-1} & { } & { } & (q^{-2}-1)q^{-(r-1)} = q^{-(r+1)}(1-q^{2}) &  (\alpha_i |\alpha_j) = -1\\ 
q^{-2} & { } & { } & (q^{-4}-1)q^{-2(r-1)} = q^{-2(r+1)}(1-q^{4}) &  (\alpha_i |\alpha_j) = -2\\ 
q^{-3} & { } & { } & (q^{-6}-1)q^{-3(r-1)} = q^{-3(r+1)}(1-q^{6}) &  (\alpha_i |\alpha_j) = -3
\end{array}\right.$$

On the algebra $\N_q^-$ let us define the following {\it ``twisted concatenation product''}: For elements $x_{im}^-, x_{jn}^- \in \N_q^-$, we define

$$ x_{im}^- \star x_{jn}^- = \begin{cases} x_{im}^-x_{jn}^- & \text{if } i+m \geq j+n \text{ or } i+m < j+n \\
 & \text{and } (\alpha_i | \alpha_j) >0\\
q^{-(\alpha_i | \alpha_j)}x_{im}^-x_{jn}^- & \text{if } j=i+1 \text{ and } n=m+1\\
& \text{or } i=j+2 \text{ and } n=m+4\\
q^{-(\alpha_i | \alpha_j)}(x_{im}^-x_{jn}^- - x_{j,n-1}^-x_{i,m+1}^-) & \text{if } i=j+1, m+1<n, j+n-1<i+m+1\\
& \text{or } i=j+2 \text{ and } n>m+4\\
q^{-(\alpha_i | \alpha_j)}(x_{im}^-x_{jn}^- - q^{(\alpha_i | \alpha_j)}x_{i,m+1}^-x_{j,n-1}^-) & \text{if } i=j+1, m+1<n, i+m+1<j+n-1\\
& \text{or } i=j+2 \text{ and } m+2<n<m+4\\
 & \text{or } j=i+2 \text{ and } m<n\\
q^{-(\alpha_i | \alpha_j)}(x_{im}^-x_{jn}^- - x_{i,m+1}^-x_{j,n-1}) & \text{if } j=i+1 \text{ and } m+1<n\\
& \text{or } j=i+2 \text{ and } n=m+1\\
x_{im}^-x_{jn}^- - x_{i,m-1}^-x_{j,n+1} & \text{if } j=i+2 \text{ and } m=n+1
\end{cases}$$

and we define the operator $\tilde{x}_{jm}^-$ over a $\star-$monomials $(x_{i_1k_1}^- \star ( \cdots \star (x_{i_{l-1}k_{l-1}}^- \star x_{i_lk_l}^-)\cdots ))$ as $\star-$left multiplication, that is

\begin{align*}\label{staronorder}  
\tilde{x}^-_{jm}(x_{i_1k_1}^- \star ( \cdots \star (x_{i_{r-1}k_{r-1}}^- \star x_{i_lk_l}^-)\cdots ))  & := x_{jm}^-\star (x_{i_1k_1}^- \star ( \cdots \star (x_{i_{l-1}k_{l-1}}^- \star x_{i_lk_l}^-)\cdots )) \\
 & = ((x^-_{jm}\star x_{i_1k_1}^- )\star (x_{i_2k_2}^- \star ( \cdots \star (x_{i_{l-1}k_{l-1}}^- \star x_{i_lk_l}^-)\cdots )). 
\end{align*}\\

From the definition of the $\star-$product, we have that it coincides with the usual product on ordered monomials, that is,  $x_{i_1k_1}^-\cdots x_{i_lk_l}^-  = x_{i_1k_1}^-\star \cdots \star x_{i_lk_l}^- $ when $i_1 + k_1 \geq \cdots \geq i_l + k_l$. Moreover, the $\star-$product is associative over ordered monomial. In particular, for any $x_{jm}^-$ and  ordered monomial $x_{i_1k_1}^-x_{i_2k_2}^-x_{i_3k_3}^-$, we have the following associativity-like property:

\begin{align*}
\tilde{x}_{jm}^-(x_{i_1k_1}^-x_{i_2k_2}^-x_{i_3k_3}^-) &= \tilde{x}_{jm}^-(x_{i_1k_1}^-(x_{i_2k_2}^-x_{i_3k_3}^-)) \\
 &= x_{jm}^-\star (x_{i_1k_1}^-\star (x_{i_2k_2}^- \star x_{i_3k_3}^-))\\
  &=(( x_{jm}^-\star x_{i_1k_1}^- ) \star (x_{i_2k_2}^- \star x_{i_3k_3}^-))\\
  &=(( x_{jm}^-\star x_{i_1k_1}^- ) \star x_{i_2k_2}^- ) \star x_{i_3k_3}^-)
\end{align*}

 and also
 
 \begin{align*}
\tilde{x}_{jm}^-(x_{i_1k_1}^-x_{i_2k_2}^-x_{i_3k_3}^-) &= \tilde{x}_{jm}^-((x_{i_1k_1}^-x_{i_2k_2}^-)x_{i_3k_3}^-)) \\
 &= x_{jm}^-\star ((x_{i_1k_1}^-\star x_{i_2k_2}^- )\star x_{i_3k_3}^-))\\
  &=(( x_{jm}^-\star (x_{i_1k_1}^-  \star x_{i_2k_2}^- ))\star x_{i_3k_3}^-)
\end{align*}\\

So, we will say that the $\star-$product is associative up to order monomials.\\ 

We also define the operator $\tilde{\Omega}_{\psi_i}(m)$ by induction on $\star-$monomials $(x_{i_1k_1}^- \star ( \cdots \star (x_{i_{l-1}k_{l-1}} \star x_{i_lk_l}^-)\cdots ))$ as follows:

$$ \tilde{\Omega}_{\psi_i}(m)(x^-_{jk}) := \delta_{ij}\delta_{-m,k}.$$

And for a general $\star-$monomial $(x_{i_1k_1}^- \star ( \cdots \star (x_{i_{l-1}k_{l-1}}^- \star x_{i_lk_l}^-)\cdots ))$ we define

$$\tilde{\Omega}_{\psi_i}(m)((x_{i_1k_1}^- \star ( \cdots \star (x_{i_{l-1}k_{l-1}}^- \star x_{i_lk_l}^-)\cdots )))  =  \delta_{ii_1}\delta_{-m,k_1}(x_{i_2k_2}^- \star ( \cdots \star (x_{i_{l-1}k_{l-1}}^- \star x_{i_lk_l}^-)\cdots ))$$
$$+ \sum_{r\geq 0} q^{p^{mk_1}_{ii_1}}g_{i,i_1,q^{-1}}(r)(x_{i_1,m_1+r}^-\star\tilde{\Omega}_{\psi_i}(m-r)(x_{i_2k_2}^- \star ( \cdots \star (x_{i_{l-1}k_{l-1}}^- \star x_{i_lk_l}^-))) $$

where $p^{mk_1}_{i_1i} := p_{i_2, \ldots, i_l}^{k_2,\ldots, k_l}$ is defined as 

$$p^{mk_1}_{i_1i} := \begin{cases} 0 & \mbox{ if } (\alpha_i | \alpha_j)>0 \\ -(\alpha_i | \alpha_j)\min\{ \ell \, |\, \tilde{\Omega}_{\psi_i}(m-\ell)((x_{i_2k_2}^- \star ( \cdots \star (x_{i_{l-1}k_{l-1}}^- \star x_{i_lk_l}^-)) = 0 \} +1 & \mbox{ if } (\alpha_i | \alpha_j)<0\\
2 & \mbox{ if } (\alpha_i | \alpha_j)=0. \end{cases}$$ 

Note that $p^{mk_1}_{i_1i}$ exists since the $\star-$monomial $(x_{i_2k_2}^- \star ( \cdots \star (x_{i_{l-1}k_{l-1}}^- \star x_{i_lk_l}^-))$ is a product of usual monomials and the formula (5.5) in \cite{CFM01} holds.\\


Recall that order monomials are the defining elements on $\tilde{M}_q(\lambda)$, so we define the operators $\tilde{x}_{jm}^-$ and $\tilde{\Omega}_{\psi_i}(m)$ on $\tilde{M}_q(\lambda)$ as follows: For a monomial $x_{i_1k_1}^-\cdots x_{i_lk_l}^- v_{\lambda}$ we have:

\begin{equation*}\label{staronorder}  \tilde{x}^-_{jm}(x_{i_1k_1}^-\cdots x_{i_rk_r}^- v_{\lambda}) = \tilde{x}^-_{jm}(x_{i_1k_1}^- \star \cdots \star x_{i_rk_r}^- v_{\lambda})  := (x^-_{jm}\star x_{i_1k_1}^- )\star \cdots \star x_{i_rk_r}^- v_{\lambda}  \end{equation*}

\begin{align*} \tilde{\Omega}_{\psi_i}(m)(x_{i_1k_1}^-\cdots x_{i_lk_l}^- v_{\lambda}) &= \tilde{\Omega}_{\psi_i}(m)(x_{i_1k_1}^-\star\cdots\star x_{i_lk_l}^- v_{\lambda}) \\
& = \delta_{ii_1}\delta_{-m,k_1}x_{i_2k_2}^-\star\cdots\star x_{i_lk_l}^- v_{\lambda}\\  
&+ \sum_{r\geq 0} q^{p^{mk_1}_{ii_1}}g_{i,i_1,q^{-1}}(r)x_{i_1,m_1+r}\star\tilde{\Omega}_{\psi_i}(m-r)(x_{i_2k_2}^-\star\cdots\star x_{i_lk_l}^-)v_{\lambda} \\
& = \delta_{ii_1}\delta_{-m,k_1}x_{i_2k_2}^-\cdots x_{i_lk_l}^- v_{\lambda}\\  
&+ \sum_{r\geq 0} q^{p^{mk_1}_{ii_1}}g_{i,i_1,q^{-1}}(r)x_{i_1,m_1+r}\star\tilde{\Omega}_{\psi_i}(m-r)(x_{i_2k_2}^-\cdots x_{i_lk_l}^-)v_{\lambda} \\
\end{align*}

\begin{rem} Here we do not emphasize on the parenthesis because $\star-$monomials coincides with usual monomials and the associativity holds for $\star-$products up to the order. \end{rem}

\begin{prop}\label{StarOrd} $\tilde{x}_{jm}^-(x_{ik}^-)$ could be written as a linear combination of ordered monomials of length 2 with coefficients in $\Z[q]$. \end{prop}

\dem It follows directly from the definition of the $\star-$product and formula (\ref{OrdX}). \findem

\begin{prop}\label{OmStarOrd} If $x_{i_1k_1}^-\cdots x_{i_lk_l}^-$ is an ordered monomial then $\tilde{\Omega}_{\psi_i}(m)(x_{i_1k_1}^-\cdots x_{i_lk_l}^-)$ is a linear combination of ordered monomials of length l-1 with coefficients in $\Z[q]$.
\end{prop}

\dem We prove by induction on $l$. If $l=1$ or $l=2$ then the statement is clear from the definition. Let us prove for clarity when $l=3$. In this case we have

\begin{multline*}
\tilde{\Omega}_{\psi_i}(m)(x_{i_1k_1}^- x_{i_2k_2}^-x_{i_3k_3}^-) \\
= \delta_{ii_1}\delta_{-m,k_1}x_{i_2k_2}^-x_{i_3k_3}^- 
+ \sum_{r\geq 0}q^{p_{ii_1}^{mk_1}}g_{ii_1,q^{-1}}(r)x_{i_1,k_1+r}^-\star \tilde{\Omega}_{\psi_i}(m-r)(x_{i_2k_2}^-x_{i_3k_3}^-)\\
= \delta_{ii_1}\delta_{-m,k_1}x_{i_2k_2}^-x_{i_3k_3}^- + \sum_{r\geq 0}q^{p_{ii_1}^{mk_1}}g_{ii_1,q^{-1}}(r)x_{i_1,k_1+r}^-\star 
 \Bigg(\delta_{ii_2}\delta_{-m+r,k_2}x_{i_3k_3}^- \\ 
 + \sum_{r'\geq 0}q^{p_{ii_2}^{m-r,k_2}}g_{ii_2,q^{-1}}(r')x_{i_2,k_2+r'}^-\star \tilde{\Omega}_{\psi_i}(m-r-r')x_{i_3k_3}^- \Bigg) \\
= \delta_{ii_1}\delta_{-m,k_1}x_{i_2k_2}^-x_{i_3k_3}^- + \sum_{r\geq 0}q^{p_{ii_1}^{mk_1}}g_{ii_1,q^{-1}}(r)\delta_{ii_2}\delta_{-m+r,k_2}x_{i_1,k_1+r}^-\star x_{i_3k_3}^-\\
+ \sum_{r\geq 0} \sum_{r'\geq 0} \delta_{ii_3}\delta_{m-r-r',k_3}q^{p_{ii_1}^{mk_1}}q^{p_{ii_2}^{m-r,k_2}}q^{p_{ii_3}^{m-r-r',k_3}}g_{ii_1,q^{-1}}(r)g_{ii_2,q^{-1}}(r') x_{i_1,k_1+r}^-\star x_{i_2,k_2+r'}^-
\end{multline*}

which by Proposition \ref{StarOrd} is a linear combination of ordered monomials of length 2 with coefficients in $\Z[q]$. 
Assume now that $\tilde{\Omega}_{\psi_i}(m)$ acting on any monomial of length $l-1$ is a linear combination of ordered monomials of length $l-2$  with coefficients in  $\Z[q]$.Then

\begin{align*}
\tilde{\Omega}_{\psi_i}(m)(x_{i_1k_1}^-\cdots x_{i_lk_l}^-) &= \delta_{ii_1}\delta_{-m,k_1}x_{i_2k_2}^-\cdots x_{i_lk_l}^- \\
&+ \sum_{r\geq 0}q^{p_{ii_1}^{mk_1}}g_{ii_1,q^{-1}}(r)x_{i_1,k_1+r}^-\star \tilde{\Omega}_{\psi_i}(m-r)x_{i_2k_2}^-\cdots x_{i_lk_l}^-
\end{align*}

The expression in the first line is clearly ordered of length $l-1$. The expression $ x_{i_1,k_1+r}^-\star \tilde{\Omega}_{\psi_i}(m-r)x_{i_2k_2}^-\cdots x_{i_lk_l}^-$ in the second line is a linear combination of ordered monomials of length $l-1$ with coefficients in $\Z[q]$ by induction and Proposition \ref{StarOrd} which completes the proof. 
\findem

\section{The bilinear form}

We are going to define a new bilinear form on order monomials of $\mc{N}_q$. For its definition we consider the unique non-degenerate bilinear form $(-,-)$ given in Proposition \ref{bilinform}.\\

Let $x_{i_1m_1}^-\cdots x_{i_km_k}^-$, $x_{j_1n_1}^-\cdots x_{j_ln_l}^-$ be two ordered monomials of $\mc{N}_q$. We define by induction the following bilinear form: 

$$ \langle x_{im}^- , x_{jn}^- \rangle = (1, \tilde{\Omega}_{\psi_i}(-m)(x_{jn}^-)) $$

$$ \langle x_{i_1m_1}^-\cdots x_{i_km_k}^- , x_{j_1n_1}^-\cdots x_{j_ln_l}^- \rangle := \langle x_{i_2m_2}^-\cdots x_{i_km_k}^-, \tilde{\Omega}_{\psi_{i_1}}(-m_1)(x_{j_1n_1}^-\cdots x_{j_ln_l}^-) \rangle $$

\vspace{2mm}

For ${\bf i} = (i_1, \ldots, i_k), {\bf m} = (m_1, \ldots, m_k)\in \Z^{k}$ the monomial $x_{i_1m_1}^-\cdots x_{i_km_k}^-$ is denoted $x_{{\bf i}{\bf m}}^-$.

\begin{lem}\label{innerprod} Let ${\bf i} = (i_1, \ldots, i_k), {\bf m} = (m_1, \ldots, m_k)\in \Z^{k}$, ${\bf j} = (j_1, \ldots, j_l), {\bf n} = (n_1, \ldots, n_l) \in \Z^l$ such that $i_1+m_1\geq \cdots \geq i_k+m_k$, $j_1+n_1 \geq \cdots \geq j_l+n_l$ and $k>l$, then $\langle x_{{\bf i}{\bf m}}^-, x_{{\bf j}{\bf n}}^- \rangle=0$. \end{lem}

\dem The proof is by induction on $l$. First of all, we have for $l=1$,

\begin{align*}
 \langle x_{{\bf i}{\bf m}}^-, x_{j_1n_1}^- \rangle &= \langle x_{i_1m_1}x_{i_2m_2}^-\cdots x_{i_km_k}^-, x_{j_1n_1} \rangle = \langle x_{i_2m_2}^-\cdots x_{i_km_k}^-, \tilde{\Omega}_{\psi_{i_1}}(-m_1)x_{j_1n_1} \rangle \\
 &= \langle x_{i_2m_2}^-\cdots x_{i_km_k}^-, \delta_{i_1j_1}\delta_{m_1n_1}\rangle = \delta_{i_1j_1} \delta_{m_1n_1} \langle x_{i_2m_2}^-\cdots x_{i_km_k}^-, 1\rangle \\
&= \delta_{i_1j_1}\delta_{m_1n_1}\langle x_{i_3m_3}^-\cdots x_{i_km_k}^-, \tilde{\Omega}_{\psi_{i_2}}(-m_2)1 \rangle  \\
&= 0
\end{align*}

Assume the result is true for ordered monomials of length $l-1$. Then for the monomial $ x_{j_1n_1}^-\cdots x_{j_ln_l}^-$ we have

\begin{align*}
 \langle x_{{\bf i}{\bf m}}^-, x_{{\bf j}{\bf n}}^- \rangle &=  \langle x_{i_1m_1}^-\cdots x_{i_km_k}^-, x_{j_1n_1}^-\cdots x_{j_ln_l}^- \rangle = \langle x_{i_2m_2}^-\cdots x_{i_km_k}^-, \tilde{\Omega}_{\psi_{i_1}}(-m_1)x_{j_1n_1}^-\cdots x_{j_ln_l}^- \rangle \\
&= \delta_{i_1j_1}\delta_{m_1n_1}\langle x_{i_2m_2}^-\cdots x_{i_km_k}^-, x_{j_2n_2}^-\cdots x_{j_ln_l}^- \rangle  \\  
&+\sum_{r\geq 0}q^{p_{i_1j_1}^{-m_1,n_1}}g_{i_1j_1,q^{-1}}(r)\langle x_{i_2m_2}^-\cdots x_{i_km_k}^-, x_{j_1,n_1+r}\star \tilde{\Omega}_{\psi_{i_1}}(-m_1-r)x_{j_2n_2}^-\cdots x_{j_ln_l}^- \rangle 
\end{align*}

By induction, the term in the second line is zero. For the third line we have that $x_{j_1,n_1+r}\star\tilde{\Omega}_{\psi_{i_1}}(-m_1-r)x_{j_2n_2}^-\cdots x_{j_ln_l}^- $ is a linear combination of ordered monomial of length $l-1$ by Proposition \ref{OmStarOrd}. So by induction 
$$ \langle x_{i_2m_2}^-\cdots x_{i_km_k}^-, x_{j_1,n_1+r}\star\Omega_{\psi_{i_1}}(-m_1-r)x_{j_2n_2}^-\cdots x_{j_ln_l}^- \rangle = 0$$
and we are done.

\findem

\begin{lem}\label{innerxint} Let ${\bf i} = (i_1, \ldots, i_k), {\bf j} = (j_1, \ldots, j_k), {\bf m} = (m_1, \ldots, m_k), {\bf n} = (n_1, \ldots, n_k)\in \Z^{k}$, and consider ordered monomials $x_{{\bf i}{\bf m}}^- = x_{i_1m_1}^-\cdots x_{i_km_k}^-$ and $x_{{\bf j}{\bf n}}^- = x_{j_1n_1}^-\cdots x_{j_kn_k}^-$. Then $$\langle x_{{\bf i}{\bf m}}^-, x_{{\bf j}{\bf n}}^-\rangle\in\Z[q]$$ 
\end{lem}

\dem The proof goes by induction on $k$. First consider the case $k=1$. Then

$$\langle x_{i_1m_1}^-,x_{j_1n_1}^-\rangle = (1, \tilde{\Omega}_{\psi_{i_1}}(-m_1)x_{j_1n_1}^- ) = \delta_{i_1j_1}\delta_{m_1n_1} \in \Z[q].$$

\vspace{0.5cm}

We will show what happens in the case $k=2$.

\begin{align*}
\langle x_{i_1m_1}^-x_{i_2m_2}^-, x_{j_1n_1}^-x_{j_2n_2}^- \rangle &= \langle x_{i_2m_2}^-, \tilde{\Omega}_{\psi_{i_1}}(-m_1)x_{j_1n_1}^-x_{j_2n_2}^- \rangle \\
& =\delta_{i_1j_1}\delta_{m_1n_1}\langle x_{i_2m_2}^-, x_{j_2n_2}^- \rangle  \\
&+ \sum_{r \geq 0}q^{p^{-m_1,n_1}_{i_1j_1}}g_{i_1j_1,q^{-1}}(r)\langle x_{i_2m_2}^-, x_{j_1,n_1+r}^-\star \tilde{\Omega}_{\psi_{i_1}}(-m_1-r)x_{j_2n_2}^-\rangle \\
&= \delta_{i_1j_1} \delta_{m_1n_1} \delta_{i_2j_2} \delta_{m_2n_2}\\
&+  \sum_{r \geq 0}q^{p^{-m_1,n_1}_{i_1j_1}}g_{i_1j_1,q^{-1}}(r) \delta_{i_1j_2}\delta_{m_1+r,n_2} \langle x_{i_2m_2}^-, x_{j_1,n_1+r}^- \rangle \\
&= \delta_{i_1j_1} \delta_{m_1n_1} \delta_{i_2j_2} \delta_{m_2n_2}\\
&+  \sum_{r \geq 0}q^{p^{-m_1,n_1}_{i_1j_1}}g_{i_1j_1,q^{-1}}(r) \delta_{i_1j_2}\delta_{m_1+r,n_2} \delta_{i_2j_1}\delta_{m_2,n_1+r} \\
&= \delta_{{\bf i}{\bf j}}\delta_{{\bf m}{\bf n}}\\
&+  \sum_{r \geq 0}q^{p^{-m_1,n_1}_{i_1j_1}}g_{i_1j_1,q^{-1}}(r) \delta_{i_1j_2} \delta_{m_1+r,n_2} \delta_{i_2j_1}\delta_{m_2,n_1+r}.
\end{align*}

Since $q^{p^{-m_1,n_1}_{i_1j_1}}g_{i_1j_1,q^{-1}}(r) \in \Z[q]$, it follows that $\langle x_{i_1m_1}^-x_{i_2m_2}^-, x_{j_1n_1}^-x_{j_2n_2}^- \rangle \in \Z[q]$.\\

Assume by induction that the result is true for monomials of length $k-1$. Then 

\begin{align*}
 \langle x_{{\bf i}{\bf m}}^-, x_{{\bf j}{\bf n}}^- \rangle &= \langle x_{i_1m_1}^-\cdots x_{i_km_k}^-, x_{j_1n_1}^-\cdots x_{j_kn_k}^- \rangle = \langle x_{i_2m_2}^-\cdots x_{i_km_k}^-, \tilde{\Omega}_{\psi_{i_1}}(-m_1)x_{j_1n_1}^-\cdots x_{k_ln_k}^- \rangle \\
&= \delta_{i_1j_1}\delta_{m_1n_1} \langle x_{i_2m_2}^-\cdots x_{i_km_k}^-, x_{j_2n_2}^-\cdots x_{j_ln_l}^- \rangle  \\  
&+\sum_{r\geq 0}q^{p_{i_1j_1}^{-m_1,n_1}}g_{i_1j_1,q^{-1}}(r) \langle x_{i_2m_2}^-\cdots x_{i_km_k}^-, x_{j_1,n_1+r}\star\Omega_{\psi_{i_1}}(-m_1-r)x_{j_2n_2}^-\cdots x_{j_kn_k}^- \rangle 
\end{align*}

by Proposition \ref{OmStarOrd} we have that $x_{j_1,n_1+r}\star\Omega_{\psi_{i_1}}(-m_1-r)x_{j_2n_2}^-\cdots x_{j_kn_k}^- $ is a linear combination of ordered monomials of length $k-1$ with coefficients in $\Z[q]$. Hence, by induction $\langle x_{i_2m_2}^-\cdots x_{i_km_k}^-, x_{j_2n_2}^-\cdots x_{j_ln_l}^- \rangle$ and $\langle x_{i_2m_2}^-\cdots x_{i_km_k}^-, x_{j_1,n_1+r}\star\Omega_{\psi_{i_1}}(-m_1-r)x_{j_2n_2}^-\cdots x_{j_kn_k}^- \rangle$ are in $\Z[q]$. Since the coefficients in each expression are also in $\Z[q]$, the result follows. 

\findem

\begin{prop}\label{innerprod04} Let ${\bf i} = (i_1, \ldots, i_k), {\bf m} = (m_1, \ldots, m_k)\in \Z^k$ and  ${\bf j} = (j_1, \ldots, j_l), {\bf n} = (n_1, \ldots, n_l)\in \Z^{l}$, and consider ordered monomials $x_{{\bf i}{\bf m}}^-, x_{{\bf j}{\bf n}}^-$ such that $\sum_{r=1}^k (i_r+m_r) = \sum_{r=1}^l (j_r+n_r)$, then $$ \langle x_{{\bf i}{\bf m}}^-, x_{{\bf j}{\bf n}}^- \rangle = (\delta_{kl}\delta_{{\bf i}{\bf j}}\delta_{{\bf m}{\bf n}} + q\Z) \mod q^2\Z[q]$$

\end{prop}

\dem The proof is by induction on the length of the monomials. First we show that if $k\neq l$ then $\langle x_{{\bf i}{\bf m}}^-, x_{{\bf j}{\bf n}}^- \rangle =0$. When 
$k>l$ it follows from Lemma \ref{innerprod}. Now assume $k < l$. Then

\begin{align*}  \langle x_{i_1m_1}^-\cdots x_{i_km_k}^-, x_{j_1n_1}^-\cdots x_{j_ln_l}^- \rangle   &=  \langle x_{i_2m_2}^-\cdots x_{i_km_k}^-, \tilde{\Omega}_{\psi_{i_1}}(-m_1)x_{j_1n_1}^-\cdots x_{j_kn_k}^- \rangle\\
  & = \langle 1 , \tilde{\Omega}_{\psi_{i_k}}(-m_k)\cdots\tilde{\Omega}_{\psi_{i_1}}(-m_1)x_{j_1n_1}^-\cdots x_{j_ln_l}^- \rangle\\
& = (1 , \tilde{\Omega}_{\psi_{i_k}}(-m_k)\cdots\tilde{\Omega}_{\psi_{i_1}}(-m_1)x_{j_1n_1}^-\cdots x_{j_ln_l}^- )\\
& = (\tilde{\Omega}_{\psi_{i_k}}(-m_k)\cdots\tilde{\Omega}_{\psi_{i_1}}(-m_1)x_{j_ln_l}^-\cdots x_{j_ln_l}^- ,1)\\
\end{align*}

The second and third equalities follow from definition. The last equality follows by the symmetry of the inner product $(-,-)$ given in Proposition \ref{bilinform}. By Proposition \ref{OmStarOrd}, we know that $\tilde{\Omega}_{\psi_{i_k}}(-m_k)\cdots\tilde{\Omega}_{\psi_{i_1}}(-m_1)x_{j_1n_1}^-\cdots x_{j_ln_l}^- $ is a linear combination of  ordered monomial of lenght $l-k>1$. So, by \cite{CFM01} Proposition 6.0.6 and the fact that $\Omega$ operators kills the element $1$, we get $(\tilde{\Omega}_{\psi_{i_k}}(-m_k)\cdots\tilde{\Omega}_{\psi_{i_1}}(-m_1)x_{j_1n_1}^-\cdots x_{j_ln_l}^- ,1)=0$.\\

So, we can now assume that $k=l$.
First, we will show the cases of monomials of length 1 and 2 to show the idea of the proof and then we will go for the general case. 

For $k=1$ we have:

\begin{align*}
\langle x_{i_1m_1}^-,x_{j_1n_1}^- \rangle & =  (1,\tilde{\Omega}_{\psi_{i_1}}(-m_1)x_{j_1n_1}^-)  \\
& = \delta_{i_1j_1}\delta_{m_1n_1}(1, 1) + \sum_{r\geq 0} q^{p_{i_1j_1}^{-m_1,n_1}}g_{i_1j_1,q^{-1}}(r)(1, x_{j_1,n_1+r}^-\star \Omega_{\psi_{i_1}}(-m_1-r)\cdot 1) \\
&  = \delta_{i_1j_1}\delta_{m_1n_1}, \
\end{align*}
since $ \Omega_{\psi_{i_1}}(-m_1-r)\cdot 1 = 0$.
Assume now $k=2$. 

\begin{align*}
& \langle x_{i_1m_1}^-x_{i_2m_2}^-, x_{j_1n_1}^-x_{j_2n_2}^- \rangle = \langle x_{i_2m_2}^-, \tilde{\Omega}_{\psi_{i_1}}(-m_1)x_{j_1n_1}^-x_{j_2n_2}^-\rangle \\
& = \delta_{i_1j_1}\delta_{m_1n_1}\langle x_{i_2m_2}^-, x_{j_2n_2}^- \rangle \\
&+ \sum_{r \geq 0}q^{p^{-m_1,n_1}_{i_1j_1}}g_{i_1j_1,q^{-1}}(r) \langle x_{i_2m_2}^-, x_{j_1,n_1+r}^-\star \tilde{\Omega}_{\psi_{i_1}}(-m_1-r)x_{j_2n_2}^- \rangle \\
& \equiv \delta_{{\bf i}{\bf j}}\delta_{{\bf m}{\bf n}} + val_{i_1j_1}^{m_1n_1} \langle x_{i_2m_2}^-, x_{j_1,n_1+\varepsilon_{i_1j_1}}^-\star \tilde{\Omega}_{i_1}(-m_1-\varepsilon_{i_1j_1})x_{j_2n_2}^- \rangle \mod q^2\Z[q] \\
& \equiv \delta_{{\bf i}{\bf j}}\delta_{{\bf m}{\bf n}} + val_{i_1j_1}^{m_1n_1} \delta_{i_1j_2}\delta_{m_1+\varepsilon_{i_1j_1},n_2}\langle x_{i_2m_2}^-, x_{j_1,n_1+\varepsilon_{i_1j_1}}^- \rangle \mod q^2\Z[q] \\
& \equiv q^{p^{-m_1,n_1}_{i_1j_1}}q^{p^{-m_2,n_2}_{i_2j_2}}\delta_{{\bf i}{\bf j}}\delta_{{\bf m}{\bf n}} + val_{i_1j_1}^{m_1n_1}\delta_{i_1j_2} \delta_{m_1+\varepsilon_{i_1j_1},n_2} \delta_{i_2,j_1} \delta_{m_2,n_1+\varepsilon_{i_1j_1}}\mod q^2\Z[q]
\end{align*}

where 

$$ \varepsilon_{i_1j_1} = \begin{cases} 1 & (\alpha_{i_1} |\alpha_{j_1})>0 \\ p_{i_1j_1}^{-m_1,n_1}  - 2& (\alpha_{i_1} |\alpha_{j_1})<0 \\ 0 & (\alpha_{i_1} |\alpha_{j_1})=0  \end{cases}   $$

and $val_{i_1j_1}^{m_1n_1} = q^{p^{-m_1,n_1}_{i_1j_1}}g_{i_1j_1,q^{-1}}(\varepsilon_{i_1j_1}) \mod q^2\Z[q]$, that is 

$$ val_{i_1j_1}^{m_1n_1} = \begin{cases} 1 & (\alpha_{i_1} |\alpha_{j_1})>0 \\ q &   (\alpha_{i_1} |\alpha_{j_1})<0 \\ 0 &  (\alpha_{i_1} |\alpha_{j_1})=0 \end{cases}$$ 

The second term is clearly zero if $(\alpha_{i_1} |\alpha_{j_1})=0$ (this parts becomes zero when we mod out by $q^2$). Let us assume that $(\alpha_{i_1} |\alpha_{j_1})\neq 0$  holds. In this case, the second term is different from zero only when: $i_1=j_2$, $m_1+\varepsilon_{i_1j_1}=n_2$, $i_2=j_1$ and $m_2=n_1+\varepsilon_{i_1j_1}$ then 

$$ i_1 + m_1 \geq i_2 + m_2 = j_1 +n_1 + \varepsilon_{i_1j_1} > j_1+n_1 \geq j_2+n_2 = i_1+m_1+\varepsilon_{i_1j_1}$$

which is not possible. Then we have

\begin{align*}
\langle x_{i_1m_1}^-x_{i_2m_2}^-, x_{j_1n_1}^-x_{j_2n_2}^- \rangle & \equiv \delta_{{\bf i}{\bf j}}\delta_{{\bf m}{\bf n}} \mod q^2\Z[q] 
\end{align*}

Let us now consider the general case.

\begin{align*} & \langle x_{i_1m_1}^-\cdots x_{i_km_k}^-, x_{j_1n_1}^-\cdots x_{j_kn_k}^- \rangle   =  \langle x_{i_2m_2}^-\cdots x_{i_km_k}^-, \tilde{\Omega}_{\psi_{i_1}}(-m_1)x_{j_1n_1}^-\cdots x_{j_kn_k}^- \rangle\\
& = \delta_{i_1j_1} \delta_{m_1n_1} \langle x_{i_2m_2}^-\cdots x_{i_km_k}^-, x_{j_2n_2}^-\cdots x_{j_kn_k}^- \rangle \\
& + \sum_{r\geq 0} q^{p^{-m_1,n_1}_{i_1j_1}} g_{i_1j_1,q^{-1}}(r) \langle x_{i_2m_2}^-\cdots x_{i_km_k}^-, x_{j_1,n_1+r}^- \star \tilde{\Omega}_{\psi_{i_1}}(-m_1-r)x_{j_2n_2}^-\cdots x_{j_kn_k}^- \rangle \\
&\equiv   \delta_{{\bf i}{\bf j}}\delta_{{\bf m}{\bf n}}  + q\Z\\
&+ val_{i_1j_1}^{m_1n_1} \langle x_{i_2m_2}^-\cdots x_{i_km_k}^-, x_{j_1,n_1+\varepsilon_{i_1j_1}}^- \star \tilde{\Omega}_{\psi_{i_1}}(-m_1-\varepsilon_{i_1j_1})x_{j_2n_2}^-\cdots x_{j_kn_k}^- \rangle \mod q^2\Z[q]
\end{align*}

Let us analyze the term $ val_{i_1j_1}^{m_1n_1} \langle x_{i_2m_2}^-\cdots x_{i_km_k}^-, x_{j_1,n_1+\varepsilon_{i_1j_1}}^- \star \tilde{\Omega}_{\psi_{i_1}}(-m_1-\varepsilon_{i_1j_1})x_{j_2n_2}^-\cdots x_{j_kn_k}^- \rangle $. If $(\alpha_{i_1} | \alpha_{j_1})= 0$ it is zero since $ val_{i_1j_1}^{m_1n_1} = 0$.  If $(\alpha_{i_1} | \alpha_{j_1})< 0$,  this term reduces to $q\Z \mod q^2\Z[q]$ by  Lemma \ref{innerxint}. Now assume $(\alpha_{i_1} | \alpha_{j_1}) > 0$.
If $i_1=j_1$ we will get that 

\begin{align*} 
& val_{i_1j_1}^{m_1n_1} \langle x_{i_2m_2}^-\cdots x_{i_km_k}^-, x_{j_1,n_1+\varepsilon_{i_1j_1}}^- \star \tilde{\Omega}_{\psi_{i_1}}(-m_1-\varepsilon_{i_1j_1})x_{j_2n_2}^-\cdots x_{j_kn_k}^- \rangle \\ &= \langle x_{i_2m_2}^-\cdots x_{i_km_k}^-, x_{j_1,n_1+1}^- \star \tilde{\Omega}_{\psi_{i_1}}(-m_1-1)x_{j_2n_2}^-\cdots x_{j_kn_k}^- \rangle\\
& \equiv \delta_{i_1j_2}\delta_{m_1+1,n_2} \langle x_{i_2m_2}^-\cdots x_{i_km_k}^-, x_{j_1,n_1+1}^- \star x_{j_3n_3}^-\cdots x_{j_kn_k}^- \rangle  \\
& + val_{i_1j_2}^{m_1+1,n_2} \langle x_{i_2m_2}^-\cdots x_{i_km_k}^-, x_{j_1,n_1+1}^- \star x_{j_2,n_2+\varepsilon_{i_1j_2}}^- \star \tilde{\Omega}_{\psi_{i_1}}(-m_1-1)x_{j_3n_3}^-\cdots x_{j_kn_k}^- \rangle \mod q^2\Z[q].
\end{align*}

Since $j_1+n_1+1\geq j_3+n_3$ we have by induction that 

\begin{align*} & \delta_{i_1j_2}\delta_{m_1+1,n_2} \langle x_{i_2m_2}^-\cdots x_{i_km_k}^-, x_{j_1,n_1+1}^- \star x_{j_3n_3}^-\cdots x_{j_kn_k}^- \rangle   \\
& = \langle x_{i_2m_2}^-\cdots x_{i_km_k}^-, x_{j_1,n_1+1}^-  x_{j_3n_3}^-\cdots x_{j_kn_k}^- \rangle \\
& \equiv  \delta_{i_1j_2}\delta_{m_1+1,n_2} \delta_{i_2j_1}\delta_{m_2,n_1+1}\delta_{i_3j_3}\delta_{m_3n_3}\cdots\delta_{i_kj_k}\delta_{m_kn_k} + q\Z \mod q^2\Z[q]
\end{align*}

and the first term is different from zero only if $i_1=j_2$, $i_2=j_1$, $m_1+1=n_2$, $m_2=n_1+1$ and $i_l=j_l$, $m_l=n_l$ for $3\leq l\leq k$. But in this case

$$ i_1+m_1\geq i_2+m_2 = j_1+n_1+1 > j_1 +n_1 \geq j_2+n_2 = i_1+m_1+1$$

which is impossible. So, when $i_1=j_1$ we get

\begin{align*} 
& val_{i_1j_1}^{m_1n_1} \langle x_{i_2m_2}^-\cdots x_{i_km_k}^-, x_{j_1,n_1+\varepsilon_{i_1j_1}}^- \star \tilde{\Omega}_{\psi_{i_1}}(-m_1-\varepsilon_{i_1j_1})x_{j_2n_2}^-\cdots x_{j_kn_k}^- \rangle \\
\equiv & val_{i_1j_2}^{m_1+1,n_2} \langle x_{i_2m_2}^-\cdots x_{i_km_k}^-, x_{j_1,n_1+1}^- \star x_{j_2,n_2+\varepsilon_{i_1j_2}}^- \star \tilde{\Omega}_{\psi_{i_1}}(-m_1-1)x_{j_3n_3}^-\cdots x_{j_kn_k}^- \rangle +q\Z  \mod q^2\Z[q]
\end{align*}

\vspace{3mm}

Now applying the same argument as above for $i_1$ and $j_2$ we get

\begin{align*} 
&  val_{i_1j_2}^{m_1+1,n_2} \langle x_{i_2m_2}^-\cdots x_{i_km_k}^-, x_{j_1,n_1+1}^- \star x_{j_2,n_2+\varepsilon_{i_1j_2}}^- \star \tilde{\Omega}_{\psi_{i_1}}(-m_1-1)x_{j_3n_3}^-\cdots x_{j_kn_k}^- \rangle \\
& \equiv   val_{i_1j_3}^{m_1+2,n_3} \langle x_{i_2m_2}^-\cdots x_{i_km_k}^-, x_{j_1,n_1+1}^- \star x_{j_2,n_2+1}^- \star x_{j_3,n_3+\varepsilon_{i_1j_3}}^-\star  \tilde{\Omega}_{\psi_{i_1}}(-m_1-1)x_{j_3n_3}^-\cdots x_{j_kn_k}^- \rangle \\
 & +q\Z \mod q^2\Z[q]
\end{align*}

\vspace{3mm}

Continuing this way, we see that the case to be analyzed is when $i_1=j_1=\cdots = j_k$.  However, this case reduces to the same proof as for $\hat{sl}_2$ given in \cite{CFM03}. \\

Hence, we see that  

$$ \langle x_{i_1m_1}^-\cdots x_{i_km_k}^-, x_{j_1n_1}^-\cdots x_{j_kn_k}^-\rangle \equiv ( \delta_{{\bf i}{\bf j}}\delta_{{\bf m}{\bf n}} +q\Z)  \mod q^2\Z[q] $$

\findem

\section{Crystal lattices}

Recall that $\mb{A}_0=\C[q^{1/2}]_{(q)}$ the ring of rational functions in $q^{1/2}$ regular at $0$. Let $\pi = \{ -\alpha + n\delta | \alpha\in \Delta_{0,+}, n\in\Z \}\cup \{0\}$. 

\begin{defi}\label{CL} Let $M$ be a $U_q(\hat{\g})$-module. We call a free $\mb{A}_0$-submodule $\mc{L}$ of $M$ an {\it imaginary crystal lattice} of $M$ if the following holds:
\begin{enumerate}
\item $\C(q^{1/2})\otimes_{\mb{A}_0}\mc{L}\cong M$.
\item $\mc{L} \cong \oplus_{\lambda\in \pi}\mc{L}_{\lambda}$ and $\mc{L}_{\lambda} = \mc{L}\cap M_{\lambda}$.
\item $\tilde{\Omega}_{\psi_i}(m)\mc{L}\subseteq \mc{L}$ and $\tilde{x}_{im}^-\mc{L}\subseteq \mc{L}$, for $i\in I_0$ and $m\in\Z$.
\end{enumerate}\end{defi}

We will show now that imaginary crystal lattices exist. Let $\lambda \in P$ such that $\lambda(c)=0$ and $\lambda(h_i)\neq 0$, $i\in I_0$, then $\tilde{M}_q(\lambda)$ is a simple reduced imaginary Verma module. \\

Consider the following $\mb{A}_0$-module:

$$ \mc{L}(\lambda) := \bigoplus_{\substack{k\geq 0 \\i_1 + m_1\geq \cdots\geq i_k + m_k\\ i_l,m_l \in\Z}} \mb{A}_0  x_{i_1 m_1}^-\cdots x_{i_k m_k}^- v_{\lambda}$$

\vspace{2mm}

Properties $(1)$ and $(2)$ clearly hold. To show that the property $(3)$ in Definition $\ref{CL}$ holds, consider the element $x_{i_1 m_1}^-\cdots x_{i_k m_k}^- v_{\lambda}$ in $\mc{L}(\lambda)$. Using induction on the length $k$ of the monomial we will show that $\tilde{\Omega}_{\psi_i}(m)x_{i_1 m_1}^-\cdots x_{i_k m_k}^- v_{\lambda} \in \mc{L}(\lambda)$. It holds for $k=1$ since\\

\begin{equation}\label{OxinB} \tilde{\Omega}_{\psi_i}(m)x_{i_1m_1}^-v_{\lambda} = \delta_{ii_1}\delta_{-m,m_1} \in \mc{L}(\lambda). \end{equation}

Assume that the result is true for monomials of length $k-1$. Then by definition we have

\begin{multline*} \tilde{\Omega}_{\psi_i}(m)(x_{i_1m_1}^-\cdots x_{i_km_k}^- v_{\lambda}) \\=  \delta_{ii_1}\delta_{-m,m_1}x_{i_2m_2}^-\cdots x_{i_km_k}^- v_{\lambda} + \sum_{r\geq 0} q^{p^{mm_1}_{ii_1}}g_{i,i_1,q^{-1}}(r)x_{i_1,m_1+r}\star \tilde{\Omega}_{\psi_i}(m)(x_{i_2m_2}^-\cdots x_{i_km_k}^-) v_{\lambda} 
 \end{multline*}

and $q^{p^{mm_1}_{i1_1}}g_{i_1i,q^{-1}}(r)  \in \mb{A}_0$ for $0\leq r\leq \varepsilon_{ii_1}$. By Proposition \ref{OmStarOrd}, $x_{i_1,m_1+r}\star \tilde{\Omega}_{\psi_i}(m)(x_{i_2m_2}^-\cdots x_{i_km_k}^-) v_{\lambda}$ is a linear combination of ordered monomials of length $k-1$ with coefficients in $\Z[q]$. So by induction we have that $\tilde{\Omega}_{\psi_i}(m)x_{i_1 m_1}^-\cdots x_{i_k m_k}^- v_{\lambda} \in \mc{L}(\lambda)$.\\

Proposition \ref{StarOrd} implies that $ \tilde{x}_{im}^-(x_{i_1 m_1}^-\cdots x_{i_k m_k}^- v_{\lambda}) \in \mc{L}(\lambda)$. 
This proves that $\mc{L}(\lambda)$ is an imaginary crystal lattice of $M$. \\

\begin{prop} $$ \mc{L}(\lambda) = \{ u\in \tilde{M}_q(\lambda) | \langle u, \mc{L}(\lambda)\rangle \subset \mb{A}_0  \}  $$ \end{prop}

\dem Let $L_1=\{ u\in \tilde{M}_q(\lambda) | \langle u, \mc{L}(\lambda) \rangle \subset \mb{A}_0  \}$. From Proposition \ref{innerprod04} it follows that $\mc{L}(\lambda)\subseteq L_1$, because the monomials in consideration are ordered. Let see the converse. Suppose $u\in L_1$, then $u\in\tilde{M}_q(\lambda)$ and $\langle u,v \rangle \in \mb{A}_0$, for any $v\in\mc{L}(\lambda)$.\\

Let $u=\sum a_{{\bf i}{\bf m}}(q^{1/2})x_{{\bf i}{\bf m}}^-v_{\lambda}$ and take the biggest $k\geq 0$ such that there exists ${\bf j}, {\bf n} \in \Z^{l}$, for some $l$, with the property that $a_{{\bf j}{\bf n}}(q^{1/2}) = q^{-k/2}b_{{\bf j}{\bf n}}(q^{1/2})$, i.e.,  $a_{{\bf j}{\bf m}}(q^{1/2})$ has a pole of order $k$. The remaining $a_{{\bf i}{\bf m}}(q^{1/2})$ has poles of order at most $k$, i.e., $a_{{\bf i}{\bf m}}(q^{1/2}) = q^{-l_{{\bf i}{\bf m}}/2}b_{{\bf i}{\bf m}}(q^{1/2})$ where $l_{{\bf i}{\bf m}} \leq k$ and $b_{{\bf i}{\bf m}}(q^{1/2})\in\C[q^{1/2}]$. Then by Proposition \ref{innerprod04} we have

$$ \langle u, x_{{\bf j}{\bf n}}^-v_{\lambda} \rangle =  q^{-k/2}b_{{\bf j}{\bf n}}(q^{1/2})(1 + qa + q^2 s_{{\bf j}{\bf n}}(q)) + \sum_{{\bf i}\neq{\bf j}, {\bf m}\neq{\bf n}}q^{-l_{\bf i}/2}b_{{\bf i}{\bf m}}(q^{1/2})(qc + q^2s_{{\bf i}{\bf m}}(q)) $$

where $s_{{\bf j}{\bf n}}(q), s_{{\bf i}{\bf m}}(q)\in\Z[q]$, $a,c\in\Z$. By hypothesis $\langle u, x_{{\bf j}{\bf m}}^-v_{\lambda} \rangle \in\mb{A}_0$, which is possible when $k=0$, so $l_{\bf i}=0$ and $b_{{\bf i}{\bf m}}q^2s_{\bf i}(q)\in\C[q^{1/2}]$. Hence, we get $a_{{\bf i}{\bf m}}\in\mb{A}_0$ and $u\in\mc{L}(\lambda)$ as desired.

%
%
%

\begin{defi}\label{deficrystalbase}
An {\it imaginary crystal basis} of a reduced imaginary Verma module $\tilde{M}_q(\lambda)$ is a pair $(\mc{L},\mc{B})$ satisfying:
\begin{enumerate}
\item $\mc{L}$ is an imaginary crystal lattice of $\tilde{M}_q(\lambda)$.
\item $\mc{B}$ is a $\C$-basis of $\mc{L}/q\mc{L} \cong \C \otimes_{\mb{A}_0}\mc{L}$.
\item $\mc{B} = \cup_{\mu\in\pi} \mc{B}_{\mu}$ where $\mc{B}_{\mu}= \mc{B}\cap(\mc{L}_{\mu}/q\mc{L}_{\mu})$.
\item $\tilde{x}_{im}^-\mc{B}\subset \pm\mc{B}\cup \{0\}$ and $\tilde{\Omega}_{\psi_i}(m)\mc{B}\subset \pm\mc{B}\cup \{0\}$, for $i \in I_0$ and $m \in \Z$.
\item For $m\in\Z$ and $i\in I_0$ if $\tilde{\Omega}_{\psi_i}(-m)b\neq 0$ and $\tilde{x}_{im}^-b\neq 0$ for $b\in\mc{B}$, then $\tilde{x}_{im}^-\tilde{\Omega}_{\psi_i}(-m)b= \tilde{\Omega}_{\psi_i}(-m)\tilde{x}_{im}^-b$.
\end{enumerate}\end{defi}

For $\lambda\in\mf{h}^*$ define

$$ \mc{B}(\lambda) = \Bigg \{ x_{i_1 m_1}^-\cdots x_{i_k m_k}^- v_{\lambda} + q\mc{L}(\lambda) \in \mc{L}(\lambda)/q\mc{L}(\lambda) \Bigg |  \begin{array}{c} i_1+ m_1 \geq \cdots \geq i_k + m_k,\\ 
m_1, \ldots, m_k \in \Z,  i_1, \ldots, i_k \in I_0 \end{array} \Bigg \}.  $$

\vspace{0.5cm}

\begin{thm}\label{CrystalbaseSimple} If $\lambda\in\mf{h}^*$ such that $\lambda(c)=0$ and $\lambda(h_i)\neq 0$ for all $i\in I_0$, then $(\mc{L}(\lambda), \mc{B}(\lambda))$ is an imaginary crystal bases for $\tilde{M}_q(\lambda)$. \end{thm}

\dem Properties $(1)$, $(2)$ and $(3)$ clearly hold. To see property $(4)$ holds, assume that $b= x_{i_1 m_1}^-\cdots x_{i_k m_k}^- v_{\lambda} + q\mc{L}(\lambda) \in \mc{B}(\lambda)$, and consider $\tilde{x}_{im}^-b$. If it is already ordered, we are done. In other cases, by the definition of the $\star-$procuct $\tilde{x}_{im}^-b + q\mc{L}(\lambda) = 0$ in the second, third and fifth case. For the fourth and sixth cases, we can order the expression by induction thanks to the $\star-$associativity over ordered monomials and Equation \ref{OrdX}. Hence, $\tilde{x}_{im}^-b+ q\mc{L}(\lambda) \in \mc{B}(\lambda)\cup \{0\}$. \\


Note that $\tilde{x}_{im}^-b = 0$ in $\mc{B}(\lambda)$ for example when there exists $j\geq 1$ such that $i=i_1=\cdots = i_j$ and $m_j=m+j$. Other possibility is when there exists $j>1$, such that $i=i_1=\cdots = i_j$, $m_l>m+l$, for $1<l<j-1$ and $m_j=m+j$.



\vspace{2mm}

Next we will see that $ \tilde{\Omega}_{\psi_i}(m)b+ q\mc{L}(\lambda)  \in \mc{B}(\lambda) \cup \{ 0 \} $ by induction on the length of the monomial $b$. Using Equation (\ref{OxinB}) we see that $ \tilde{\Omega}_{\psi_i}(m)x_{i_1m_1}^- + q\mc{L}(\lambda)  \in \mc{B}(\lambda) \cup \{ 0 \} $.  Consider now

\begin{align*} & \tilde{\Omega}_{\psi_i}(m)x_{i_1 m_1}^-\cdots x_{i_k m_k}^- v_{\lambda} = \delta_{ii_1}\delta_{-m,m_1}x_{i_2m_2}^-\cdots x_{i_km_k}^- v_{\lambda} \\ 
& +  \sum_{r \geq 0} q^{p^{mm_1}_{i1_1}}g_{i_1i,q^{-1}}(r) x_{i_1,m_1+r}^- \star \tilde{\Omega}_{\psi_i}(m-r) x_{i_2 m_2}^-\cdots x_{i_k m_k}^-v_{\lambda}  \\ 
& \equiv  \delta_{ii_1}\delta_{m,-m_1}x_{i_2m_2}^-\cdots x_{i_km_k}^- v_{\lambda}  + wal_{ii_1}^{mm_1} x_{i_1, m_1+\varepsilon_{ii_1}}^-\star \tilde{\Omega}_{\psi_i}(m-\varepsilon_{ii_1})x_{i_2m_2}^-\cdots x_{i_km_k}^- v_{\lambda} \mod q\mc{L}(\lambda) \\
\end{align*}

where $\varepsilon_{ii_1}$ is defined as in Proposition \ref{innerprod04} and 

$$ wal_{ii_1}^{m_1n_1} = \begin{cases} -1 & (\alpha_{i_1} | \alpha_{j_1}) >0\\ 0 &(\alpha_{i_1} | \alpha_{j_1}) \leq 0 \end{cases}$$

Then, if $(\alpha_{i_1} | \alpha_{j_1}) \leq 0$ we get $\tilde{\Omega}_{\psi_i}(m)x_{i_1 m_1}^-\cdots x_{i_k m_k}^- v_{\lambda} = 0 \in \mc{B}(\lambda)\cup \{0\}$. In other case, 

\begin{align*} & \tilde{\Omega}_{\psi_i}(m)x_{i_1 m_1}^-\cdots x_{i_k m_k}^- v_{\lambda} \\ 
& \equiv  \delta_{m,-m_1}x_{i_2m_2}^-\cdots x_{i_km_k}^- v_{\lambda}  - x_{i_1, m_1+1}^-\star \tilde{\Omega}_{\psi_i}(m-1)x_{i_2m_2}^-\cdots x_{i_km_k}^- v_{\lambda} \mod q\mc{L}(\lambda). 
\end{align*}

By induction $ \tilde{\Omega}_{\psi_i}(m-1)x_{i_2m_2}^-\cdots x_{i_km_k}^- v_{\lambda} \mod q\mc{L}(\lambda) \neq 0$ only if $i=i_2=\cdots=i_k$ and in this situation we get

\begin{align*} & \tilde{\Omega}_{\psi_i}(m)x_{i_1 m_1}^-\cdots x_{i_k m_k}^- v_{\lambda} \\ 
& \equiv  \delta_{m,-m_1}x_{i_2m_2}^-\cdots x_{i_km_k}^- v_{\lambda}  + x_{i_1, m_1+1}^-\star \tilde{\Omega}_{\psi_i}(m-1)x_{i_2m_2}^-\cdots x_{i_km_k}^- v_{\lambda} \mod q\mc{L}(\lambda)  \\
&\equiv  \sum_{j=1}^{k+1} (-1)^{j-1} \delta_{m-j+1,-m_j} x_{i_1, m_1+1}^-\star x_{i_2,m_2+1}^-\star \cdots \star x_{i_{j-1},m_{j-1}+1}^-\star \cdots \star x_{i_km_k+1}^- \mod q\mc{L}(\lambda),
\end{align*}

where $x_{0,m_0+1} =1$.
Note that $x_{i_1, m_1+1}^-\star x_{i_2,m_2+1}^-\star \cdots \star x_{i_{j-1},m_{j-1}+1}^-\star \cdots \star x_{i_km_k+1}^-$ is already order because $x_{i_1 m_1}^-\cdots x_{i_k m_k}^-$ is ordered. So, $\tilde{\Omega}_{\psi_i}(m)x_{i_1 m_1}^-\cdots x_{i_k m_k}^- v_{\lambda} \mod q\mc{L}(\lambda) \in \mc{B}(\lambda)\cup\{0\}$. \\

For $(5)$ consider $b= x_{i_1 m_1}^-\cdots x_{i_k m_k}^- v_{\lambda} + q\mc{L}(\lambda) \in \mc{B}(\lambda)$ assume that $\tilde{\Omega}_{\psi_i}(-m)b\neq 0$ and $\tilde{x}_{im}b\neq 0$. As we saw, this is just possible if $i=i_1=\cdots=i_k$ and there does not exist $j$ such that $m_j=m+j$. In such a case it reduces the proof to the $\hat{{sl}}_2$ situation which is already proved in \cite{CFM03}. Hence we are done. 

\findem

\section{The category $\mc{O}^q_{red,im}$}

Let $G_q$ be the quantized Heisenberg subalgebra generated by $h_{in}$ and $\gamma$, for $i\in I_0$, $n\in \Z \setminus \{0\}$. \\

We will say that a $U_q(\hat{\g})$-module V is $G_q$-compatible if:\\

\begin{enumerate}
    \item [(i)] $V$ has a decomposition $V=T(V)\oplus TF(V)$ where $T(V)$ and $TF(V)$ are non-zero $G_q$-modules, called, respectively, torsion and torsion free module associated to $V$.
    \item [(ii)] $h_{im}$ for $i\in I_0$, $m\in \Z\setminus\{0\}$ acts bijectively on $TF(V)$, i.e., they are bijections on $TF(V)$.
    \item [(iii)] $TF(V)$ has no non-zero $U_q(\hat{\g})$-submodules.
    \item [(iv)] $G_q\cdot T(V)=0$.
\end{enumerate}

\vspace{3mm}

Consider the set 

$$\mf{h}^*_{q,red} = \{ \lambda \in P \; | \; \lambda(c)=0, \lambda(h_i)\neq 0, i\in I_0 \}$$

The category $\mc{O}^q_{red,im}$ is defined as the category whose objects are $U_q(\hat{\g})$-modules $M$ such that \\

\begin{enumerate}
    \item $M$ is $\hat{\h}^*_{q,red}$-diagonalizable, that means,  $$ M = \bigoplus_{\nu \in \hat{\h}^*_{q,red}} M_{\nu}, \mbox{  where  } M_{\nu} = \{ m\in M | K_im=q^{\lambda(h_i)}m, Dm = q^{\lambda(d)}m, i\in I_0 \}$$
    \item For any $i\in I_0$ and any $n\in \Z$, $x^+_{in}$ acts locally nilpotently.
    \item $M$ is $G_q$-compatible.
    \item the morphisms in $\mc{O}^q_{red,im}$ are $U_q(\hat{\g})$-homomorphisms.
\end{enumerate}

\vspace{3mm}

Reduced imaginary Verma modules belongs to $\mc{O}^q_{red,im}$. Indeed, for $\tilde{M}_q(\lambda)$ consider $T(\tilde{M}_q(\lambda)) = \C(q^{1/2}) v_{\lambda}$ and $TF(V) = \bigoplus_{k\in\Z, n_1, \ldots, n_N\in\Z_{\geq0}} \tilde{M}_q(\lambda)_{\lambda+k\delta - n_1\alpha_1 - \ldots -n_N\alpha_N}$, and at least one $n_j\neq 0$. Moreover, direct sums of reduced imaginary Verma modules belongs to $\mc{O}^q_{red,im}$. \\

Due to \cite{FutGM}, we can deform  $U(\hat{\g})$-imaginary Verma modules preserving weight space decompositions and weight multiplicities. So the proof of the following theorem is completely analogues to the one presented in Theorem 5.1, Theorem 5.3 and Proposition 5.4 in \cite{AFO} which we sketch here. 

\begin{thm}
\begin{itemize}
\item[(1)]
If $\lambda,\mu \in \hat{\h}^*_{q,red}$ then $\Ext_{\mc{O}_{red,im}^q}^1(\tilde{M}_q(\lambda), \tilde{M}_q(\mu)) = 0$.
\item[(2)] If $M$ is an irreducible module in the category $\mc{O}_{red,im}^q$, then $M\cong \tilde{M}_q(\lambda)$ for some $\lambda \in \hat{\mf{h}}_{q,red}^*$. Moreover, if $N$ is an arbitrary object of $\mc{O}^q_{red,im}$ then $N\cong \bigoplus_{\lambda_i \in \hat{\mf{h}}^*_{q,red}} \tilde{M}(\lambda_i)$, for some $\lambda_i's$.
\end{itemize}
\end{thm}

\dem Suppose there exists an extension $M$ that fits in the following short exact sequence: \[ \xymatrix{0 \ar[r] & \tilde{M}_q(\lambda) \ar[r]^{\iota} & M \ar[r]^{\pi} & \tilde{M}_q(\mu) \ar[r] & 0 }. \]
If $\lambda$ and $\mu$ just differ by a multiple of the null root, we will get two vectors in $M$, which are annihilated  by $x_{in}^{+}$ for $i\in I_0$ and $n\in \Z$ and, by the $G_q$-compatibility, they are isolated. Hence, they are highest weight vectors, which generates two irreducible submodules isomorphic to $\tilde{M}_q(\lambda)$ and $\tilde{M}_q(\mu)$ and then the extension splits. \\

Assume now $\mu = \lambda +k\delta- \sum_{i=1}^N s_i\alpha_i $, for $s_i\in\Z$, $k\in \Z$ and with not all $s_i$ equal to zero. First of all, let $s_i\in\Z_{\geq 0}$ and let $\overline{v}_\mu$ be a preimage under $\pi$ of a highest weight vector $v_{\mu}$ of $\tilde{M}_q(\mu)$ of highest weight $q^\mu$. Because $x_{in}^+ v_{\mu} = G_q v_\mu = 0$, for any $i\in I_0$ and $n\in \Z$, it is possible to show that $T(M) = \C(q^{1/2})v_\mu$. Moreover, for $i\in I_0$ and $s\in \Z\setminus\{0\}$ we have that $\pi(h_{is}\overline{v}_{\mu}) = h_{is}v_\mu =0$, then $h_{is}\overline{v}_{\mu} \in \tilde{M}_q(\lambda)$, but it is just possible that $h_{is}\overline{v}_{\mu}=0$ and so $G_q\overline{v}_{\mu}=0$.\\

Since the operators $x_{in}^+$ are locally nilpotent, we can see that $x_{in}^+\overline{v}_\mu=0$ for all $i \in I_0$ and $n\in \Z$. Indeed, If this is not the case, there exists $j\in I_0$ and $m\in \Z$ such that $0\neq x_{in}^+\overline{v}_\mu \in \tilde{M}_q(\lambda)$. If we fixed $j$ we have a $U_q(\mf{sl}_2)$-subalgebra $U_q(j)$ such that the $U_q(j)$-submodule of $M$ generated by $\overline{v}_\mu$, say $M_q(j)$, it is an extension of reduced imaginary $U_q(j)$-Verma modules, one of them of weight $q^\mu$. So, because $M\in \mc{O}_{red,im}^q$ we have that $M_q(j)\in \mc{O}_{red,im}^q(U_q(j))$, but this is a semisimple category by \cite{CFM03}, hence $x_{in}^+\overline{v}_\mu=0$ for all $i$ and $n$. Therefore, $\overline{v}_{\mu}$ generates a $U_q(\g)$-submodule of $M$ isomorphic to $\tilde{M}_q(\mu)$ and the short exact sequence splits. \\

In case $s_i\in\Z_{\leq 0}$ for all $i$ and at least one different from $0$, because $\tilde{M}_q(\mu)$ is irreducible and $\tilde{M}(\lambda)$ is a $U_q(\g)$-submodule of $M$, the short exact sequence splits completing the proof of statement (1).\\

Assume now that $M \in \mc{O}^q_{red,im}$ is irreducible. Let $v\in T(M)$ be a nonzero element of weight $\lambda\in\mf{h}^*_{q,red}$. For each $i\in I_0$ let $p_i \in \Z_{>0}$ be the minimum possible integer such that $(x^+_{i0})^{p_i}v=0$. If all $p_i=1$ we see that $x^+_{in}v=0$ for all $i\in I_0$ and $n\in\Z\setminus \{0\}$. Hence, we have an epimorphism $\tilde{M}_q(\lambda) \twoheadrightarrow M$, so $M\cong \tilde{M}(\lambda)$.\\

On the other hand, assume there exists at least one $p_i$ such that $p_i>1$. Then, there exists $\ell\in \Z_{> 0}$ and a nonzero element $w_{\bf i}$, for ${\bf i} = (i_1, \ldots, i_{\ell})\in I_0^{\ell}$ such that $x^+_{jn}w_{\bf i}=0$, for any $j\in I_0$ and $n\in \Z$. Consider the $G_q$-submodule  $W_{\bf i} = U(G_q)w_{\bf i}$ of $M$ and the induced module $I(W_{\bf i}) = \Ind_{\mc{B}_q}^{U_q(\hat{\g})} W_{\bf i}$, where $\mc{B}_q$ is generated by $x_{in}^+$, $K_i^{\pm 1}$, $D^{\pm 1}$, $h_{im}$ and $\gamma$ for $i\in I_0$, $m,n\in\Z$, $m\neq 0$. Since $M$ is irreducible, it is a quotient of $I(W_{\bf i})$. If $w_{\bf i}\in T(M)$, we have $W_{\bf i} = \C(q^{1/2}) w_{\bf i}$, and so $M$ is a quotient of $I(W_{\bf i}) = \tilde{M}_q(\lambda)$ and we are done. In case $w_{\bf i}\notin T(M)$, we get a contradiction. We conclude that $M\cong \tilde{M}_q(\lambda)$ for some $\lambda \in \hat{\mf{h}}_{q,red}^*$.\\

Let $N$ be an arbitrary module in $\mc{O}_{red,im}^q$ and $v\in T(N)$ is nonzero. Let $w_{\bf i}$ and $W_{\bf i}$ as above. Then we have two possibilities: either $w_{\bf i}\notin T(N)$ or $w_{\bf i} \in T(N)$. In the first case, we get a proper $G_q$-submodule of $TF(N)$ which is not possible. In the second case, $W_{\bf i} = \C(q^{1/2}) w_{\bf i} \subseteq T(N)$ and for some $\lambda_{\bf i}$, $I(W_{\bf i}) \cong \tilde{M}_q(\lambda_{\bf i})$ is a $U_q(\hat{\g})$-submodule of $N$. Then, any non-zero element of $T(N)$ generates an irreducible reduced imaginary Verma module which is a $U_q(\hat{\g})$-submodule of $N$ and because there are no extensions between them, they are direct summands of $N$.\\
\findem

For $i\in I_0$ and $n,m \in \Z$ we have defined operators $\tilde{x}_{in}$ and $\tilde{\Omega}_{\psi_i}(m)$ on irreducible reduced imaginary Verma modules. So, due to the above theorem we have the following:

\begin{thm}
The operators $\tilde{x}_{in}$ and $\tilde{\Omega}_{\psi_i}(m)$ are well defined on objects in the category $\mc{O}^q_{red,im}$.
\end{thm}
\findem

For an object $M\in \mc{O}^q_{red,im}$, we can define imaginary crystal lattices and bases analogous to the Definitions \ref{CL} and \ref{deficrystalbase}. In what follows we will prove that crystal basis of irreducible reduced imaginary Verma modules extend to direct sums of these modules. We also prove a partial converse of this statement. The proofs are analogues to the proofs in the $\hat{sl}_2$-case given in \cite{CFM04}.\\

Now suppose $M\in \mc{O}^q_{red,im}$, then there exists $\lambda_k\in \hat{\mf{h}}^*_{q,red}$ for $k\in J$ ($J$ an index set) such that $M\cong \bigoplus_{k \in J} \tilde{M}(\lambda_k)$. For $k\in J$, let $(\mc{L}(\lambda_k), \mc{B}(\lambda_k))$ be the imaginary crystal basis of $\tilde{M}_q(\lambda_k)$ given in Theorem \ref{CrystalbaseSimple}. Set $\mc{L} = \bigoplus_{k\in J} \mc{L}(\lambda_k)$ and $\mc{B} = \bigsqcup_{k\in J}\mc{B}(\lambda_k)$. 

\begin{thm}\label{directcrystal} Let $M\in \mc{O}^q_{red,im}$ such that $M\cong \bigoplus_{\lambda_k \in J} \tilde{M}(\lambda_k)$ as above. Then the pair $(\mc{L},\mc{B})$ is an imaginary crystal basis for $M$.
\end{thm}

\dem We need to check that the five conditions in Definition \ref{deficrystalbase} hold. For the first one we need to see that $\mc{L}$ is an imaginary crystal lattice. Clearly, $\C(q^{1/2})\otimes_{\mb{A}_0}\mc{L} \cong \bigoplus_{k\in J} \C(q^{1/2})\otimes_{\mb{A}_0}\mc{L}(\lambda_k) \cong  \bigoplus_{k\in J} \tilde{M}_q(\lambda_k) \cong M$. The third property  follows directly because $\tilde{\Omega}_{\psi_i}(m)(\mc{L}) = \tilde{\Omega}_{\psi_i}(m)(\bigoplus_{k\in J} \mc{L}(\lambda_k)) = (\bigoplus_{k\in J} \tilde{\Omega}_{\psi_i}(m)\mc{L}(\lambda_k)) \subseteq \bigoplus_{k\in J} \mc{L}(\lambda_k) = \mc{L}$ and similarly we have that $\tilde{x}_{im}^-(\mc{L}) = \tilde{x}_{im}^-(\bigoplus_{k\in J} \mc{L}(\lambda_k)) = (\bigoplus_{k\in J} \tilde{x}_{im}^-\mc{L}(\lambda_k)) \subseteq  \bigoplus_{k\in J} \mc{L}(\lambda_k) = \mc{L}$, for any $i\in I_0$ and $m\in \Z$. So let us see the second property in Definition \ref{CL}. For this we first show that $\mc{L}_{\mu} = (\bigoplus_{k\in J} \mc{L}(\lambda_k))_{\mu} = \bigoplus_{k\in J} \mc{L}(\lambda_k)_{\mu}$, where $\mc{L}(\lambda_k)_{\mu} = \mc{L}(\lambda_k) \cap M_{\mu}$. \\

In fact, assume $u\in (\bigoplus_{k\in J} \mc{L}(\lambda_k))_{\mu}$, then $u=\sum_{k\in J} u_k$ where $u_k\in \mc{L}(\lambda_k)$ and $K_iu = q^{\mu(h_i)}u$ for any $i\in I_0$. Since $\mc{L}(\lambda_k) = \bigoplus_{\mu_k \in \pi} \mc{L}(\lambda_k)_{\mu_k}$, we can write $u_k=\sum_{\mu_k\in \pi} u_{k,\mu_k}$ where $u_{k,\mu_k}\in \mc{L}(\lambda_k)_{\mu_k}$ and $K_iu_{k,\mu_k} = q^{\mu_k(h_i)}u_{k,\mu_k}$ for any $i\in I_0$. Hence, we have the following

\begin{align*}
K_iu &= \sum_{k\in J} K_iu_k = \sum_{k\in J} \sum_{\mu_k \in \pi} K_iu_{k,\mu_k} =  \sum_{k\in J} \sum_{\mu_k \in \pi} q^{\mu_k(h_i)}u_{k,\mu_k}\\
\parallel & \\
q^{\mu(h_i)}u &= \sum_{k\in J} q^{\mu(h_i)}u_k = \sum_{k\in J} \sum_{\mu_k \in \pi} q^{\mu(h_i)}u_{k,\mu_k}.
\end{align*}

Hence $\sum_{k\in J} (\sum_{\mu_k \in \pi} (q^{\mu(h_i)} - q^{\mu_k(h_i)})u_{k,\mu_k})=0$. Since the sum $\bigoplus_{k\in J}\mc{L}(\lambda_k)$ is direct we get $\mu(h_i) = \mu_k(h_k)$ for any $i\in I_0$ and $k\in J$. Then $u_{k,\mu_k}\in \mc{L}(\lambda_k)\cap M_{\mu} = \mc{L}(\lambda_k)_{\mu}$. Finally note that $\mc{L} = \bigoplus_{k\in J}\mc{L}(\lambda_k) = \bigoplus_{k\in J}\bigoplus_{\mu\in \pi} \mc{L}(\lambda_k)_{\mu} = \bigoplus_{\mu\in \pi}\bigoplus_{k\in J} \mc{L}(\lambda_k)_{\mu} = \bigoplus_{\mu\in \pi} \mc{L}_{\mu}$, where $\mc{L}_{\mu} = \mc{L}\cap M_{\mu}$. So $\mc{L}$ is an imaginary crystal lattice of $M$.\\

Let us now look at the second property in Definition \ref{deficrystalbase}. We know that $\mc{B}(\lambda_k)$ is a $\C$-basis of $\mc{L}(\lambda_k)/q\mc{L}(\lambda_k) \cong \C\otimes_{\mb{A}_0}\mc{L}(\lambda_k)$ for all $k\in J$. Now, 

\begin{align*}
\mc{L}/q\mc{L} \cong& \bigoplus_{k\in J}\mc{L}(\lambda_k)/ q\bigoplus_{k\in J}\mc{L}(\lambda_k)
\cong \bigoplus_{k\in J} (\mc{L}(\lambda_k)/q\mc{L}(\lambda_k))\\
\cong& \bigoplus_{k\in J}(\C\otimes_{\mb{A}_0}\mc{L}(\lambda_k)) \cong \C\otimes_{\mb{A}_0}\mc{L}.
\end{align*}

Hence, $\mc{L}/q\mc{L}$ has the $\C$-basis $\bigsqcup_{k\in J} \mc{B}(\lambda_k) = \mc{B}$. For the third property we have 

$$ \mc{B} = \bigsqcup_{k\in J} \mc{B}(\lambda_k) = \bigsqcup_{k\in J}\bigsqcup_{\mu\in\pi}\mc{B}(\lambda_k)_{\mu} = \bigsqcup_{\mu\in\pi}\bigsqcup_{k\in J}\mc{B}(\lambda_k)_{\mu} = \bigsqcup_{\mu\in\pi}\mc{B}_\mu$$

where $\mc{B}_{\mu} = \bigsqcup_{k\in J}\mc{B}(\lambda_k)_{\mu}$ and so

\begin{align*}
\mc{B}_{\mu} =& \bigsqcup_{k\in J}(\mc{B}(\lambda_k)\cap (\mc{L}(\lambda_k)_{\mu}/q\mc{L}(\lambda_k)_{\mu}))\\
=& (\bigsqcup_{k\in J}(\mc{B}(\lambda_k))\cap (\bigoplus_{k'\in J}\mc{L}(\lambda_{k'})_{\mu}/q\mc{L}(\lambda_{k'})_{\mu}))\\
=& \mc{B}\cap (\bigoplus_{k'\in J}\mc{L}(\lambda_{k'})_{\mu}/q\mc{L}(\lambda_{k'})_{\mu}))\\
=& \mc{B} \cap \mc{L}_{\mu}/q\mc{L}_{\mu}
\end{align*}

The fourth property holds since $\tilde{\Omega}_{\psi_i}(m)(\mc{B}) = \tilde{\Omega}_{\psi_i}(m)(\bigsqcup_{k\in J} \mc{B}(\lambda_k)) = (\bigsqcup_{k\in J} \tilde{\Omega}_{\psi_i}(m)\mc{B}(\lambda_k)) \subseteq \bigsqcup_{k\in J} \pm\mc{B}(\lambda_k)\cup\{0\} = \pm\mc{B}\cup\{0\}$ and similarly $\tilde{x}_{im}^-(\mc{B}) \subseteq \pm\mc{B}\cup\{0\}$, for any $i\in I_0$ and $m\in \Z$. \\

Finally, for the fifth property suppose $b\in \mc{B}$ such that $\tilde{\Omega}_{\psi_i}(-m)b\neq 0$ and $\tilde{x}_{im}^-b\neq 0$ for $m\in\Z$ and $i\in I_0$. Since by definition $\mc{B} = \bigsqcup_{k\in J}\mc{B}(\lambda_k)$ is a disjoint union there is unique index $k'\in J$ such that $b\in \mc{B}(\lambda_{k'})$. Since $(\mc{L}(\lambda_{k'}), \mc{B}(\lambda_{k'}))$ is an imaginary crystal basis we have that $\tilde{x}_{im}^-\tilde{\Omega}_{\psi_i}(-m)b= \tilde{\Omega}_{\psi_i}(-m)\tilde{x}_{im}^-b$ as desired. 

\findem

To prove a partial converse of the Theorem \ref{directcrystal} below, we need the following lemma.

\begin{lem}\label{lemcom} For any $i\in I_0$ and $m\in \Z$ the operators $\tilde{\Omega}_{\psi_i}(m)$ and $\tilde{x}^-_{im}$ commutes with any $U_q(\hat{\g})$-homomorphism in the category $\mc{O}^q_{red,im}$.
\end{lem}

\dem Let $M\in \mc{O}^q_{red,im}$. It is enough to prove the statement for an $U_q(\hat{\g})$-homomorphism $\varphi: \tilde{M}_q(\lambda) \to M$ for some $\lambda\in \hat{\h}^*_{q,red}$. Let $x_{i_1n_1}^-\cdots x_{i_kn_k}^-v_{\lambda}$ be an ordered monomial basis element of $\tilde{M}_q(\lambda)$, then

\begin{align*}
\varphi(\tilde{x}_{im}(x_{i_1n_1}^-\cdots x_{i_kn_k}^-v_{\lambda})) =& \varphi((x_{im}\star x_{i_1n_1}^-)\cdots \star x_{i_kn_k}^-v_{\lambda})\\
=& ((x_{im}\star x_{i_1n_1}^-)\cdots \star x_{i_kn_k}^-)\varphi(v_{\lambda})\\
=& \tilde{x}_{im}(x_{i_1n_1}^-\cdots x_{i_kn_k}^-\varphi(v_{\lambda}))\\
=& \tilde{x}_{im}(\varphi(x_{i_1n_1}^-\cdots x_{i_kn_k}^-v_{\lambda})).
\end{align*}

\begin{align*}
\varphi(\tilde{\Omega}_{\psi_i}(m)(x_{i_1n_1}^-\cdots x_{i_kn_k}^-v_{\lambda})) =& \varphi(\delta_{ii_1}\delta_{-m,n_1}x_{i_2n_2}^-\cdots x_{i_kn_k}^-v_{\lambda})\\
+& \sum_{r\geq 0}q^{p^{mn_1}_{ii_1}}g_{ii_1,q^{-1}}(r)\varphi(x_{i_1,n_1+r}^-\star\tilde{\Omega}_{\psi_1}(m-r)(x_{i_2n_2}^-\cdots x_{i_kn_k}^-v_{\lambda}))\\
=& \delta_{ii_1}\delta_{-m,n_1}x_{i_2n_2}^-\cdots x_{i_kn_k}^-\varphi(v_{\lambda})\\
+& \sum_{r\geq 0}q^{p^{mn_1}_{ii_1}}g_{ii_1,q^{-1}}(r)x_{i_1,n_1+r}^-\star\tilde{\Omega}_{\psi_1}(m-r)(x_{i_2n_2}^-\cdots x_{i_kn_k}^-\varphi(v_{\lambda}))\\
=& \tilde{\Omega}_{\psi_i}(m)(x_{i_1n_1}^-x_{i_2n_2}^-\cdots x_{i_kn_k}^-\varphi(v_{\lambda}))\\
=& \tilde{\Omega}_{\psi_i}(m)(\varphi(x_{i_1n_1}^-x_{i_2n_2}^-\cdots x_{i_kn_k}^-v_{\lambda})),
\end{align*}
which proves the statement.
\findem

\begin{thm} Let $M = M_1\oplus M_2$ where $M_1$ and $M_2$ are modules in the category
$\mc{O}^q_{red,im}$ and suppose $(\mc{L}, \mc{B})$ is an imaginary crystal basis for $M$. Furthermore, suppose that there exists $\mb{A}_0$-submodules $\mc{L}_j \subset M_j$, and subsets $\mc{B}_j \subset \mc{L}_j/q\mc{L}_j$, for $j = 1,2$ such that $\mc{L} = \mc{L}_1\oplus \mc{L}_2$ and $\mc{B} = \mc{B}_1 \sqcup \mc{B}_2$. Then $(\mc{L}_j, \mc{B}_j)$ is an imaginary crystal basis of $M_j$, for $j = 1,2$.
\end{thm}

\dem It is straightforward to see that $\C(q^{1/2})\otimes_{\mb{A}_0}(\mc{L}_j)_{\mu} = (M_j)_{\mu}$, $\mu\in\pi$, $\mc{L}_j = \mc{L}\cap M_j$ and $\mc{B}_j=\mc{B}\cap(\mc{L}_j/q\mc{L}_j)$ (see Theorem 4.2.10 (2) from \cite{HK}).\\

Let see that $\mc{L}_{\mu} = (\mc{L}_1)_{\mu} \oplus (\mc{L}_2)_{\mu}$, for $\mu\in\pi$. The $``\supseteq"$ part is obvious.  For the other inclusion assume $u\in \mc{L}_\mu$, then $u=u_1+u_2$ where $u_1\in \mc{L}_1$ and $u_2\in \mc{L}_2$. Let $i\in I_0$ then $K_iu = K_iu_1 + K_iu_2 = q^{\mu(h_i)}u_1 + q^{\mu(h_i)}u_2$ and so $K_iu_1 - q^{\mu(h_i)}u_1 = -K_iu_2 + q^{\mu(h_i)}u_2 \in \mc{L}_1\cap\mc{L}_2=\{0\}$. Hence, $u_1\in(\mc{L}_1)_{\mu}$ and $u_2\in(\mc{L}_2)_{\mu}$ and we are done.\\

For $u\in \mc{L}_j\subset \mc{L} = \bigoplus_{\mu\in\pi}\mc{L}_{\mu}$ we have $u=\sum_{\mu\in\pi} u_{\mu}$ where $u_\mu\in \mc{L}_{\mu}$. We have decomposition $u_\mu = (u_1)_{\mu} + (u_2)_{\mu}$ with $(u_j)_\mu\in (\mc{L}_j)_{\mu}$, $j = 1,2$. Consequently $u-\sum_{\mu\in\pi}(u_j)_{\mu} = \sum_{\mu\in\pi}(u_k)_{\mu}$ for $j\neq k$, which implies that  $u-\sum_{\mu\in\pi}(u_j)_{\mu} \in \mc{L}_j\cap\mc{L}_k=\{0\}$. So we have $u\in \bigoplus_{\mu\in\pi}(\mc{L}_j)_{\mu}$.\\

Let $pr_j:M\to M_j$ be the canonical projection into the $j$-component. Let $u_j\in \mc{L}_j$ and recall that for any $i\in I_0$ and $m\in\Z$, $\tilde{\Omega}_{\psi_i}(m)\mc{L}\subset\mc{L}$ and $\tilde{x}^-_{im}\mc{L}\subset\mc{L}$. Then $\tilde{\Omega}_{\psi_i}(m)u_j = \overline{u}_1 + \overline{u}_2$ and $\tilde{x}^-_{im}u_j = \check{u}_1 + \check{u}_2$, where $\overline{u}_j,\check{u}_j\in\mc{L}_j$. By Lemma \ref{lemcom} we have $ \tilde{\Omega}_{\psi_i}(m)(u_1) =  \tilde{\Omega}_{\psi_i}(m)(p_1(u_1)) = p_1(\tilde{\Omega}_{\psi_i}(m)(u_1)) = \overline{u}_1$ and $ \tilde{x}^-_{im}(u_1) = \tilde{x}^-_{im}(p_1(u_1)) = p_1(\tilde{x}^-_{im}(u_1)) = \check{u}_1$. Hence $\tilde{\Omega}_{\psi_i}(m)(u_1), \tilde{x}_{im}^-(u_1) \in \mc{L}_1$. Similarly, $\tilde{\Omega}_{\psi_i}(m)(u_2), \tilde{x}_{im}^-(u_2) \in \mc{L}_2$. This shows that $\mc{L}_j$ is a crystal lattice for $j=1,2$.\\

Notice that 

$$ \mc{L}/q\mc{L} \cong \C\otimes_{\mb{A}_0} \mc{L} \cong (\C\otimes_{\mb{A}_0} \mc{L}_1) \oplus (\C\otimes_{\mb{A}_0} \mc{L}_2) \cong  (\mc{L}_1/q\mc{L}_1) \oplus (\mc{L}_2/q\mc{L}_2).$$

\vspace{0.3cm}

Using this isomorphism we have that $\mc{B}_j = \mc{B}\cap (\mc{L}_j/q\mc{L}_j)$ is a $\C$-basis of $\mc{L}_j/q\mc{L}_j \cong \C\otimes_{\mb{A}_0}\mc{L}_j$. We also have that $\mc{B} = \mc{B}_1 \sqcup \mc{B}_2$ and thus $\mc{B}_j = \bigsqcup_{\mu\in\pi} (\mc{B}_j)_{\mu}$ where $(\mc{B}_j)_{\mu} = \mc{B}\cap ((\mc{L}_j)_\mu/q(\mc{L}_j)_{\mu})$.\\

The operators $\tilde{\Omega}_{\psi_i}(m)$ and $\tilde{x}^-_{im}$, for $i\in I_0$, $m\in\Z$ leaves stable the bases $\mc{B}_j$ because of Lemma \ref{lemcom}, i.e., $\tilde{\Omega}_{\psi_i}(m)\mc{B}_j\subset\mc{B}_j\cup\{0\}$ and $\tilde{x}^-_{im}\mc{B}_j\subset\mc{B}_j\cup\{0\}$, for $j=1,2$. Finally, if $b\in \mc{B}_j$ is such that $\tilde{\Omega}_{\psi_i}(-m)b\neq 0$ and $\tilde{x}_{im}^-b\neq 0$, then $\tilde{x}_{im}^-\tilde{\Omega}_{\psi_i}(-m)b= \tilde{\Omega}_{\psi_i}(-m)\tilde{x}_{im}^-b$ since $b\in \mc{B}$.  \\

This completes the proof that $(\mc{L}_j, \mc{B}_j)$ is an imaginary crystal basis for $M_j$, $j=1,2$.

\findem

\bigskip
\begin{center}
ACKNOWLEDGEMENT
\end{center}
\smallskip
JCA is support by the FAPESP Grant 2021/13022-9. He acknowledges the hospitality of the faculty and staff in the Department of Mathematics at North Carolina State University during his visit where part of this work was done.
KCM is partially supported by Simons Foundation grant \#  636482.

\end{document}